\documentclass[psamsfonts,11pt,reqno]{amsart}
\usepackage{amsmath}
\usepackage{amsfonts}
\usepackage{amssymb,mathrsfs}

\usepackage[left=1in,right=1in,top=1in,bottom=1in]{geometry}

\usepackage[all,cmtip]{xy}

\def\N{\mathbb{N}}
\def\ap{absolutely productive}

\newtheorem{theorem}{Theorem}[section]
\newtheorem{corollary}[theorem]{Corollary}
\newtheorem{question}[theorem]{Question}
\newtheorem{proposition}[theorem]{Proposition}
\newtheorem{lemma}[theorem]{Lemma}
\newtheorem{example}[theorem]{Example}
\newtheorem{remark}[theorem]{Remark}

\newtheorem{problem}[theorem]{Problem}
\newtheorem{definition}[theorem]{Definition}
\newcommand{\Cp}[2]{C_p(#1,#2)}
\newcommand{\eq}[1]{\stackrel{#1}{\sim}}

\newcommand{\IANS}{TAP}
\newcommand{\tor}[1]{$tor(#1)$}

\newcommand{\Chom}{\mathrm{Chom}}
\renewcommand{\hat}{\widehat}

\begin{document}
\title[Group-valued continuous functions]{Group-valued continuous functions with the topology of pointwise convergence}
\renewcommand{\thefootnote}{\fnsymbol{footnote}}

\author[D. Shakhmatov]{Dmitri Shakhmatov}
\address[Dmitri Shakhmatov]{Division of Mathematics, Physics and Earth Sciences,
Graduate School of Science and Engineering\\
Ehime University, Matsuyama 790-8577, Japan}
\email{dmitri@dpc.ehime-u.ac.jp}

\author[J. Sp\v{e}v\'{a}k]{Jan Sp\v{e}v\'{a}k}
\address[Jan Sp\v{e}v\'{a}k]{}
\email{prednosta.stanice@quick.cz}


\dedicatory{Dedicated to Professor Tsugunori Nogura on the occasion of his 60th anniversary}

\thanks{The first
author was partially supported by the Grant-in-Aid for Scientific Research~(C)
no.~19540092 by the Japan Society for the Promotion of Science (JSPS)}

\thanks{This manuscript has been written as
a part of the second author's Ph. D. study at the Graduate School of
Science and Engineering of Ehime University (Matsuyama, Japan) from
October 1, 2006 till September 30, 2009. His stay at Ehime
University during this period has been generously supported by the
Ph. D. Scholarship of the Government of Japan.}

\begin{abstract}
Let $G$ be a topological group with the identity element $e$. Given
a space $X$, we denote by $\Cp{X}{G}$ the group of all continuous
functions from $X$ to $G$ endowed with the topology of pointwise
convergence, and we say that $X$ is: (a) {\em $G$-regular\/} if, for
each closed set $F\subseteq X$ and every point $x\in X\setminus F$,
there exist $f\in\Cp{X}{G}$ and $g\in G\setminus\{e\}$ such that
$f(x)=g$ and $f(F)\subseteq\{e\}$; (b) {\em $G^\star$-regular\/}
provided that there exists $g\in G\setminus\{e\}$ such that, for
each closed set $F\subseteq X$ and every point $x\in X\setminus F$,
one can find $f\in\Cp{X}{G}$ with $f(x)=g$ and $f(F)\subseteq\{e\}$.
 Spaces $X$ and $Y$ are {\em $G$-equivalent\/} provided that the
topological groups $\Cp{X}{G}$ and $\Cp{Y}{G}$ are topologically
isomorphic.

We investigate which topological properties are preserved by
$G$-equivalence, with a special emphasis being placed on
characterizing topological properties of $X$ in terms of those of
$\Cp{X}{G}$. Since $\mathbb{R}$-equivalence coincides with
$l$-equivalence, this line of research ``includes'' major topics of
the classical $C_p$-theory
 of Arhangel'ski\u{\i} as a particular case (when $G=\mathbb{R}$).

We introduce a new class of \IANS\ groups that contains all groups
having no small subgroups (NSS groups). We prove that: (i) for a
given NSS group $G$, a $G$-regular space $X$ is pseudocompact if and
only if $\Cp{X}{G}$ is \IANS, and (ii) for a metrizable NSS group
$G$, a
 $G^\star$-regular space $X$ is compact if and only if $\Cp{X}{G}$ is
a \IANS\ group of countable tightness. In particular, a Tychonoff
space $X$ is pseudocompact (compact) if and only if
$\Cp{X}{\mathbb{R}}$ is a \IANS\ group (of countable tightness).
Demonstrating the limits of the result in (i), we give an example of
a precompact TAP group $G$ and a $G$-regular countably compact space
$X$ such that $\Cp{X}{G}$ is not \IANS.

We show that Tychonoff spaces $X$ and $Y$ are
$\mathbb{T}$-equivalent if and only if their free precompact Abelian
groups are topologically isomorphic, where $\mathbb{T}$ stays for
the quotient group $\mathbb{R}/\mathbb{Z}$. As a corollary, we
obtain that $\mathbb{T}$-equivalence implies $G$-equivalence for
every Abelian precompact group $G$. We
establish that $\mathbb{T}$-equivalence preserves the following
topological properties: compactness, pseudocompactness,
$\sigma$-compactness, the property of being a Lindel\"of
$\Sigma$-space, the property of being a compact metrizable space,
the (finite) number of connected components, connectedness, total
disconnectedness. An example of $\mathbb{R}$-equivalent (that is,
$l$-equivalent) spaces that are not $\mathbb{T}$-equivalent is
constructed.
\end{abstract}

\maketitle

\medskip

 In notation and terminology we follow
\cite{Dikran} and \cite{Engelking} if not stated otherwise. {\em All
topological spaces are assumed to be Tychonoff (that is, completely
regular $T_1$ spaces)
and non-empty,
 and all topological groups are
assumed to be Hausdorff.}

 By $\mathbb{N}$ we denote the set of all
natural numbers,
$\omega$ stays for the least nonzero limit ordinal,
 $\mathbb{Z}$ is the discrete additive group of
integers, $\mathbb{R}$ is the additive group of reals with its usual
topology, $\mathbb{T}$ stays for the quotient group
$\mathbb{R}/\mathbb{Z}$, and $\mathbb{Z}(n)$ denotes the cyclic
group of order $n$ (with the discrete topology).
The identity element of a group
$G$ is denoted by $e_G$, or simply by $e$ when there is no danger of confusion.

If $G$ is a topological group, then the symbol
$\widehat{G}$ stays for the completion of $G$ with respect to
 the two-sided uniformity. If
 $G=\widehat{G}$, then $G$ is called {\em complete\/}.
  It is well-known that $\widehat{G}$ always exists,
  $\widehat{G}$ is a topological group, $G$ is dense in $\widehat{G}$, and  if
  $G$ is a dense subgroup of a complete group $H$, then $\widehat{G}=H$.
If $G$ is a subgroup of some compact group, then $G$ is called {\em
 precompact}.

Recall that a space $X$ is called {\em pathwise connected\/}
provided that for every pair of points $x,y\in X$ there exists a
continuous map $\varphi:[0,1]\to X$ from the unit interval $[0,1]$
to $X$ such that $\varphi(0)=x$ and $\varphi(1)=y$. (The image
$\varphi([0,1])$ is called a {\em path\/} between $x$ and $y$.)

\section{Introduction}
\label{introduction}

\begin{definition}
\label{function:space:definition}
{\rm
Let $X$ be a space and $G$ a topological group.
\begin{itemize}
\item[(i)]
We shall use $C(X,G)$ to denote the group of all continuous
functions from $X$ to $G$, equipped with the ``pointwise group
operations''. That is, the product of $f\in C(X,G)$ and $g\in
C(X,G)$ is the function $fg\in C(X,G)$ defined by $fg(x)=f(x)g(x)$
for all $x\in X$, and the inverse
element
of $f$ is the function
$h \in C(X,G)$
defined by
$h(x)=\left(f(x)\right)^{-1}$
for all $x\in
X$.
\item[(ii)]
The family
$$
\{W(x,U):
x\in X,\
U
\mbox{ is an open subset of }
G
\},
$$
where
$$
W(x,U)=\{f\in C(X,G): f(x)\in U\},
$$
forms a subbase of the {\em topology of pointwise convergence\/} on
$C(X,G)$. We use
the
symbol $\Cp{X}{G}$ to denote the set $C(X,G)$ endowed with this
topology.
\end{itemize}
}
\end{definition}

One can easily see that $\Cp{X}{G}$ is a topological group.

\begin{definition}
{\rm
Let $G$ and $H$ be topological groups.
\begin{itemize}
\item[(i)]
Recall that $G$ and $H$ are
said to be {\em topologically isomorphic\/} if there exists  a
bijection $f:G\to H$ which is both a group homomorphism and a
homeomorphism. We write $G\cong H$ whenever $G$ and $H$ are
topologically isomorphic.
\item[(ii)]
We say that spaces $X$ and $Y$ are
{\em $G$-equivalent\/}, and denote this by $X\eq{G}Y$, provided that
$\Cp{X}{G}\cong\Cp{Y}{G}$.
\item[(iii)]
Let
$\mathscr{C}$ be a  class of spaces.
We say that a topological property $\mathscr{E}$ is {\em preserved
by $G$-equivalence within the class $\mathscr{C}$\/} provided that
the following condition holds:
if $X\in\mathscr{C}$, $Y\in \mathscr{C}$, $X\eq{G}Y$ and $X$ has the
property $\mathscr{E}$, then $Y$ must have the property
$\mathscr{E}$ as well. The sentence ``$\mathscr{E}$ is {\em
preserved by $G$-equivalence\/}'' is used as an abbreviation for
``$\mathscr{E}$ is preserved by $G$-equivalence within the class of
Tychonoff spaces''.
\item[(iv)]
Given
a class $\mathscr{C}$ of spaces,
we say that {\em
$G$-equivalence implies $H$-equivalence within the
class $\mathscr{C}$\/}
provided that the following statement holds:
If $X\in \mathscr{C}$, $Y\in\mathscr{C}$
and
$X\eq{G}Y$, then $X\eq{H}Y$. The sentence
``$G$-equivalence implies $H$-equivalence'' shall be used as an abbreviation for
``$G$-equivalence implies $H$-equivalence within the class of Tychonoff spaces''.
\end{itemize}
}
\end{definition}

In \cite{Markov} Markov has introduced
the
free
topological group $F(X)$
of
a space $X$
and defined
spaces $X$ and $Y$
to be
{\em $M$-equivalent\/} if $F(X)\cong F(Y)$.
Thereafter,
a significant effort
went into an investigation of how topological properties of $F(X)$
depend on those of $X$, as well as which topological properties are
preserved by $M$-equivalence.

Every continuous function $f:X\to G$
from a space $X$ to a topological group $G$
can be (uniquely) extended to a continuous group homomorphism
$\hat{f}: F(X)\to G$.
This elementary fact (with $\mathbb{T}$ as $G$) was applied by Graev
to show that the closed unit interval and the circle are not
$M$-equivalent \cite{Graev}.
Tkachuk noticed in \cite{Tkachuk} that
$M$-equivalence
implies $G$-equivalence for every Abelian topological group $G$.
He then applied this observation to $G=\mathbb{Z}(2)$ to show that
connectedness
is preserved by $M$-equivalence \cite{Tkachuk}.

Later on, many properties of $M$-equivalence were discovered by
means of
the notion of $l$-equivalence; see \cite{A, AT}. Recall that
spaces $X$ and $Y$ are called {\em $l$-equivalent\/}
provided that $\Cp{X}{\mathbb{R}}$ and $\Cp{Y}{\mathbb{R}}$ are {\em
topologically isomorphic as topological vector spaces\/}. A
fundamental observation pertinent to the subject of this paper has
been made in \cite{Tkachuk} by Tkachuk:
spaces $X$ and $Y$
are $l$-equivalent if and only if $C_p(X,\mathbb{R})$ and
$C_p(Y,\mathbb{R})$ are {\em topologically isomorphic as topological
groups\/}. In other words, $l$-equivalence of
spaces
coincides with their $\mathbb{R}$-equivalence (in our notation). A
far reaching conclusion that one might get from this fact is that,
despite a significant emphasis on the {\em topological vector
space\/} structure commonly placed in the $C_p$-theory \cite{A},
this structure is largely irrelevant to the study of the notion of
$l$-equivalence, and in fact may as well be replaced by the {\em
topological group\/} structure. It is this conclusion that led us to
an idea of introducing the general notion of $G$-equivalence, for
an arbitrary topological group $G$.

This opens up a
topic of studying the properties of
the topological group $\Cp{X}{G}$, for a given space $X$ and a
topological group $G$. Let us outline major problems that appear to
be of particular interest in this new area of research.

\begin{problem}
\label{problem:1}
Given a topological group $G$, characterize topological properties of
a space $X$ in terms of algebraic and/or topological properties of
$\Cp{X}{G}$.
\end{problem}
\begin{problem}
\label{problem:2}
Given a topological group $G$, a class $\mathscr{C}$ of spaces and a
topological property $\mathscr{E}$, investigate when the property $\mathscr{E}$
is preserved
by $G$-equivalence within the class $\mathscr{C}$.
\end{problem}

In the particular case when $G=\mathbb{R}$, these two problems are well-known
(and major) problems of the $C_p$-theory. Therefore, one can view Problems
\ref{problem:1} and \ref{problem:2} as a natural generalization of
major
topics of the $C_p$-theory to the case of an arbitrary topological group $G$.

Since the $C_p$-theory provides a large supply of topological properties
that are preserved by $\mathbb{R}$-equivalence (that is,
$l$-equivalence)
within the class of Tychonoff spaces,
the following
particular version of Problem \ref{problem:2} seems to be worth
studying:
\begin{problem}
Let $G$ be one of the ``important'' topological groups such as, for
example, the circle group $\mathbb{T}$, the dual group
$\mathbb{Q}^*$
 of the discrete group $\mathbb{Q}$ of
rational numbers, or the group $\mathbb{Z}_p$ of $p$-adic integers.
Assume also that $\mathscr{E}$ is a topological property preserved
by $\mathbb{R}$-equivalence within (some subclass of) the class of
Tychonoff spaces. Is it then true that $\mathscr{E}$ is also
preserved by $G$-equivalence
within
(an  appropriate subclass of) the
class of Tychonoff spaces?
\end{problem}

One can
formulate the most ambitious version of Problem \ref{problem:2}:
\begin{problem}
\label{problem:3}
Let  $\mathscr{C}$  be
a class of spaces and $\mathscr{E}$ a topological property.
Describe the class $\mathbf{G}(\mathscr{C},\mathscr{E})$ of topological groups $G$ such that the property $\mathscr{E}$
is preserved
by $G$-equivalence within the class $\mathscr{C}$.
\end{problem}

As may be expected, Problem \ref{problem:3} turns out to be difficult even in
the case of major topological properties $\mathscr{E}$ such as compactness and pseudocompactness.

\begin{problem}
Given a class $\mathscr{C}$ of spaces and topological groups $G$,
$H$, when does $G$-equivalence imply $H$-equivalence within the
class $\mathscr{C}$?
\end{problem}

In this manuscript we build a foundation for studying these
problems.

Section \ref{basic:results} collects necessary preliminaries and basic results,
most of which are elementary and have counterparts in the classical $C_p$-theory.
Example \ref{negative:example} demonstrates that, in order to obtain nice results, it is important
to have enough continuous maps from a space $X$ to a given topological group $G$. This fact naturally leads to an
introduction of three
notions of ``regularity'', $G$-regularity, $G^\star$-regularity and
$G^{\star\star}$-regularity,
 in Definition \ref{def:G:regular}. One of the basic results
(Proposition \ref{product:equivalence}) states that $G$-equivalence
and $H$-equivalence combined together imply $(G\times
H)$-equivalence. The converse implication fails in general (Example
\ref{GxH:does:not:imply:G}). The main goal of
Section~\ref{section:3} is to show that $(G\times H)$-equivalence
does imply both $G$-equivalence and $H$-equivalence for
``sufficiently different'' topological groups $G$ and $H$ (Theorems
\ref{thm:4.4} and \ref{thm:4.5}).

In Section \ref{section:4} we introduce a new class of topological groups
(that we call \IANS\ groups) and prove that every group without small subgroups (an NSS group) is \IANS; see Theorem \ref{NSSjeTAP}. The class of \IANS\ groups has many common properties with that of NSS groups.
For example, this class is closed under
taking
subgroups and finite products, and  a \IANS\ group does not
contain any subgroup topologically isomorphic to an infinite product of non-trivial topological groups
(Proposition \ref{h}). Every topological group without non-trivial convergent sequences is \IANS\
(Proposition \ref{psc:TAP:not:NSS}), so the class of \IANS\ groups contains many ``peculiar'' topological groups.

In Section \ref{section:5} we show that a space $X$ is pseudocompact if and only if $\Cp{X}{\mathbb{R}}$ is a
\IANS\ group (Theorem \ref{Rpseudo}), thereby providing
a short, elementary, ``group-theoretic'' proof
of the result of Arhangel'ski\u{\i} about preservation of pseudocompactness by $l$-equivalence.

In Section \ref{pseu}
we further generalize Theorem \ref{Rpseudo} by proving that, for every NSS group $G$, a  $G$-regular space $X$
is pseudocompact if and only if $\Cp{X}{G}$ is \IANS; see Theorem \ref{main}.
Emphasizing the limits of this result,
we construct a precompact \IANS\ group $G$
and
 a countably compact $G^\star$-regular space $X$ such that $\Cp{X}{G}$ is not
\IANS\ (Theorem \ref{separately:continuous:example}).

The main result of Section \ref{seccomp} is Theorem
\ref{char:of:compactness} saying that, for a metrizable NSS group
$G$, a $G^\star$-regular space $X$ is compact if and only if
$\Cp{X}{G}$ is a \IANS\ group of countable tightness. In particular,
$G$-equivalence  preserves compactness within the class of
$G^\star$-regular spaces for every NSS metric group $G$ (Corollary
\ref{nademnou}). The classes of topological groups $G$ such that
$G$-equivalence preserves compactness and pseudocompactness,
respectively, are closed under taking finite powers (Corollary
\ref{productOK}) but are not closed under taking finite products
(Example \ref{GxH:does:not:imply:G}). Moreover, we give an example
demonstrating that $G$-equivalence can preserve both compactness and
pseudocompactness without $G$ being NSS, or even \IANS\ (Example
\ref{E2}).

Section \ref{secdis} provides some sufficient conditions on a topological group $G$
that guarantee that $G$-equivalence preserves total disconnectedness, connectedness and (finite) number of connected components.

In Section \ref{scr:G:eq} we recall a general categorical machinery that
leads to a
definition of a free object $F_\mathscr{G}(X)$ of a space $X$ in a given class $\mathscr{G}$ of
topological groups (closed under
taking
products and subgroups), and we define spaces $X$ and $Y$ to be  $\mathscr{G}$-equivalent
provided that $F_\mathscr{G}(X)\cong F_\mathscr{G}(Y)$. When $\mathscr{G}$ is the class of all topological groups
(all topological Abelian groups, respectively), the free object $F_\mathscr{G}(X)$ coincides with the free
topological group (respectively, the free Abelian topological group) of a space $X$ in the sense of
Markov, and $\mathscr{G}$-equivalence coincides with the classical $M$-equivalence ($A$-equivalence, respectively).
When $G\in\mathscr{G}$ is Abelian, then  $\mathscr{G}$-equivalence implies $G$-equivalence (Corollary \ref{A_G implikuje G}).

Section \ref{T:section} is devoted to the study of properties of $\mathbb{T}$-equivalence.
A non-trivial connection with the previous section is based on the so-called ``precompact
duality theorem'' (Theorem \ref{duality}) that allows us to prove that $\mathbb{T}$-equivalence
coincides with $\mathscr{P}$-equivalence, where $\mathscr{P}$ is the class of all precompact Abelian
groups (Corollary \ref{PjeT}). Combining this with Corollary \ref{A_G implikuje G},
we conclude
that $\mathbb{T}$-equivalence
implies $G$-equivalence for
{\em every\/} precompact Abelian group $G$ (Corollary \ref{T:implies:any:precompact}).
Theorem \ref{properties:preserved:by:T} lists major topological properties
that are preserved by $\mathbb{T}$-equivalence. As a consequence, all these properties are
also preserved by $\mathscr{G}$-equivalence whenever $\mathbb{T}\in \mathscr{G}$
(Corollary \ref{scrG:equivalance:when:scrG:contains:T}).
In particular, it follows that
total disconnectedness is preserved by $A$-equivalence and $M$-equivalence
(Corollary \ref{total:disconnectedness:M:A}),
which seems to be a new result.
Since $l$-equivalence (aka $\mathbb{R}$-equivalence) does not preserve connectedness, while
$\mathbb{T}$-equivalence does, it follows that $l$-equivalence does not imply $\mathbb{T}$-equivalence (Proposition \ref{R:does:not:imply:T}).

Section \ref{open:questions} lists some concrete open problems that are related to our
results.

\section{Basic results}
\label{basic:results}

\begin{example}
\label{negative:example} {\rm Let $G$ be a topological group with
the trivial connected component (for example, a zero-dimensional
group). Then $\Cp{X}{G}\cong
G$ for every connected space $X$. In particular, any two connected
spaces $X$ and $Y$ are
$G$-equivalent, and so most major topological properties are not
preserved by $G$-equivalence. }
\end{example}

This example clearly demonstrates
that, in order to obtain meaningful theorems about preservation of some topological
property by $G$-equivalence within the class $\mathscr{C}$, one has to require each
member of $\mathscr{C}$ to have ``sufficiently many'' continuous
functions to the target topological group $G$. Our next definition exhibits three possible ways of doing so.

\begin{definition}
\label{def:G:regular}
{\rm Given a topological group $G$,
we
say that a space $X$ is:
\begin{enumerate}
\item[(i)]{\em$G$-regular} if, for
each closed set $F\subseteq X$
and every point $x\in X\setminus F$,
there
exist
$f\in\Cp{X}{G}$
and $g\in G\setminus\{e\}$ such that $f(x)=g$ and
$f(F)\subseteq\{e\}$;
\item[(ii)]{\em $G^\star$-regular} if there exists $g\in G\setminus\{e\}$ such
that, for every closed set $F\subseteq X$ and each point $x\in X\setminus F$, one
can find $f\in\Cp{X}{G}$ such that $f(x)=g$ and
$f(F)\subseteq\{e\}$;
\item[(iii)]
{\em $G^{\star\star}$-regular} provided that, whenever
$F$ is a closed subset of $X$,
$x\in X\setminus F$
and $g\in G$,
 there exists $f\in\Cp{X}{G}$
such that $f(x)=g$ and $f(F)\subseteq\{e\}$.
\end{enumerate}
}
\end{definition}

It is clear that
\begin{equation}
\label{eq:G:regularity:implications}
X \mbox { is } G^{\star\star}\mbox{-regular}
\to
X \mbox { is } G^{\star}\mbox{-regular}
\to
X \mbox { is } G\mbox{-regular.}
\end{equation}

Since the topological group $G=\mathbb{R}\times\mathbb{Z}(2)$ is not connected,
 and a
continuous image of a connected space is connected,
one can easily see that
that unit interval $[0,1]$ is $G^\star$-regular but not
$G^{\star\star}$-regular,
so the first implication in \eqref{eq:G:regularity:implications}
cannot be reversed. The authors have no example of a group $G$
witnessing that the second implication in
\eqref{eq:G:regularity:implications} cannot be reversed, although
they are convinced that such an example must exist;
see Question \ref{question:G:regularity}.

Our next proposition describes three obvious cases when some kind of
$G$-regularity ``comes for free'':
\begin{proposition}\label{Gregualrity}
Let $X$ be a space and $G$ a topological group.
\begin{itemize}
\item[(i)]
If $G$ is pathwise connected, then $X$ is  $G^{\star\star}$-regular.
\item[(ii)]
If $G$ contains a homeomorphic copy of the unit interval $[0,1]$,
then $X$ is $G^\star$-regular.
\item[(iii)] If $X$ is zero-dimensional in the sense of $\mathrm{ind}$, then
$X$ is $G^{\star\star}$-regular.
\end{itemize}
In particular, in all three cases, $X$ is $G$-regular by
\eqref{eq:G:regularity:implications}.
\end{proposition}

It should be noted that our terminology differs from that of
\cite{hern}, where a pair $(X,G)$ consisting of a space $X$ and a
topological group $G$ is called $G$-regular if it satisfies the
condition (iii) of Definition \ref{def:G:regular}. The same
manuscript \cite{hern} states explicitly (but using different
terminology) item (i) of Proposition \ref{Gregualrity}.

For a cardinal $\tau\ge 1$ we denote by $D_\tau$ the discrete space of size $\tau$.

\begin{proposition}
\label{proposition:4.4}
\label{group.of.countable.power}
Let $G$ be a topological group and $\tau\ge 1$ be a cardinal.
\begin{itemize}
\item[(i)]
Spaces $X$ and $Y$ are
$G^\tau$-equivalent if and only if $X\times D_\tau$ and $Y\times D_\tau$ are
$G$-equivalent.
\item[(ii)]
If $\tau$ is infinite, then every space $X$ is $G^\tau$-equivalent to $X\times D_\tau$.
\end{itemize}
\end{proposition}
\begin{proof}
(i) follows from
$C_p(X\times D_\tau,G)\cong C_p\left(X,G^{D_\tau}\right)\cong C_p(X, G^\tau)$
and $C_p(Y\times D_\tau,G)\cong C_p\left(Y,G^{D_\tau}\right)\cong C_p(Y, G^\tau)$.

(ii) follows from
$C_p\left(X\times D_\tau, G^\tau\right)\cong
C_p\left(X, (G^\tau)^{D_\tau}\right)\cong C_p\left(X, G^{\tau\times
D_\tau}\right)\cong C_p(X, G^\tau)$.
\end{proof}

Applying Propositions \ref{Gregualrity}(iii) and
\ref{group.of.countable.power}(ii), we get the following

\begin{corollary}
\label{Gomegadoesntpreserve}
\label{notcon}
For every topological group $G$,
a singleton and the countable discrete space $D_\omega$ are
$G^{\star\star}$-regular and $G^\omega$-equivalent.
In particular,
$G^\omega$-equivalence
preserves neither
pseudocompactness (compactness),
nor (pathwise) connectedness
within the class of
 $G^{\star\star}$-regular spaces.
\end{corollary}

\begin{corollary}
\label{connectedness:and:finite:powers}
If $G$-equivalence preserves the finite number of connected components,
then so does $G^k$-equivalence for all $k\in\mathbb{N}\setminus\{0\}$.
\end{corollary}
\begin{proof}
Fix $k\in\mathbb{N}\setminus\{0\}$,
 and suppose that $X$ and $Y$ and $G^k$-equivalent.
Then $X\times D_k$ and $Y\times D_k$ are $G$-equivalent by
Proposition \ref{group.of.countable.power}(i). By the assumption of
our corollary, $X\times D_k$ and $Y\times D_k$ have the same (finite) number
of connected components. Clearly, this implies that
 $X$ and $Y$ must also have the
same number of connected components.
\end{proof}

\begin{proposition}
\label{product:equivalence}
Let $\{G_i:i\in I\}$ be a family of topological groups.
If spaces $X$ and $Y$ are $G_i$-equivalent
for all $i\in I$, then $X$ and $Y$ are also $\left(\prod_{i\in I}
G_i\right)$-equivalent.
\end{proposition}
\begin{proof}
$C_p\left(X,\prod_{i\in I} G_i\right)\cong \prod_{i\in I}
C_p(X,G_i)$.
\end{proof}

\begin{corollary}
\label{comparing:powers}
Let $G$ be a topological group.
\begin{itemize}
\item[(i)]
$G$-equivalence implies $G^\kappa$-equivalence for every cardinal
$\kappa\ge1$.
\item[(ii)]
Suppose that $\tau$ and $\kappa$ are cardinals such that
$1\le\tau\le\kappa$ and $\kappa\ge \omega$.
Then $G^\tau$-equivalence implies $G^\kappa$-equivalence.
\end{itemize}
\end{corollary}

\begin{proposition}
\label{closed:copy:of:G:in:Cp}
$\Cp{X}{G}$ contains a closed subgroup topologically isomorphic to $G$.
\end{proposition}

\begin{definition}
{\rm
For a topological group $G$, a topological property $\mathscr{E}$ and a class $\mathscr{C}$ of spaces define
\begin{align*}
\mathbf{I}(G,\mathscr{E}, \mathscr{C})=\{\tau:\  & \tau\ge 1
\mbox{ is a cardinal such that }
\\
&
G^\tau
\mbox{-equivalence preserves property }
\mathscr{E}
\mbox{ within the class }
\mathscr{C}\}.
\end{align*}
}
\end{definition}

\begin{lemma}
\label{non:empty:implies:1}
Let $G$ be a topological group, $\mathscr{E}$ a topological property
and $\mathscr{C}$ a class of spaces.
Then
$\mathbf{I}(G,\mathscr{E}, \mathscr{C})\not=\emptyset$ implies
$1\in \mathbf{I}(G,\mathscr{E}, \mathscr{C})$.
\end{lemma}
\begin{proof}
Suppose $\mathbf{I}(G,\mathscr{E}, \mathscr{C})\not=\emptyset$, and
let $\kappa\in \mathbf{I}(G,\mathscr{E}, \mathscr{C})$.
Assume
that $X\in\mathscr{C}$ and $Y\in\mathscr{C}$ are $G$-equivalent spaces such that $X$ has property $\mathscr{E}$.
Then $X$ and $Y$
are also $G^\kappa$-equivalent by
Corollary \ref{comparing:powers}(i).
Since $\kappa\in  \mathbf{I}(G,\mathscr{E}, \mathscr{C})$ and $X$ has property $\mathscr{E}$,
$Y$ must have it as well. Thus, $1\in \mathbf{I}(G,\mathscr{E},\mathscr{C})$.
\end{proof}

\begin{lemma}
\label{1:implies:tau} Assume that  $G$ is a topological group,
$\tau\ge 1$ is a cardinal, $\mathscr{E}$ is a topological property
and $\mathscr{C}$ is a class of spaces satisfying the following
conditions:
\begin{itemize}
\item[(i)]
a space
$X$ has
property $\mathscr{E}$ if and only if the product $X\times D_\tau$ has property $\mathscr{E}$;
\item[(ii)]
if $X\in \mathscr{C}$, then $X\times D_\tau\in \mathscr{C}$.
\end{itemize}
Then $1\in \mathbf{I}(G,\mathscr{E},\mathscr{C})$ implies $\tau\in \mathbf{I}(G,\mathscr{E},\mathscr{C})$.
\end{lemma}
\begin{proof}
Assume that $X\in\mathscr{C}$ and $Y\in\mathscr{C}$ are $G^\tau$-equivalent spaces such that
$X$ has property $\mathscr{E}$.
Then $X\times D_\tau$ and $Y\times D_\tau$ are
$G$-equivalent by Proposition \ref{proposition:4.4}(i).
Furthermore, $X\times D_\tau\in\mathscr{C}$ and $Y\times D_\tau\in\mathscr{C}$
by (ii).
Since $X\times D_\tau$ has property $\mathscr{E}$ by (i),
and
$1\in \mathbf{I}(G,\mathscr{E},\mathscr{C})$, we conclude that $Y\times D_\tau$ must have property $\mathscr{E}$.
Applying (i) once again, we conclude that $Y$ has property $\mathscr{E}$.
This proves that $\tau\in \mathbf{I}(G,\mathscr{E},\mathscr{C})$.
\end{proof}

\begin{proposition}
Let $G$ be a topological group, $\mathscr{E}$ a topological property and
$\mathscr{C}$ a class of spaces satisfying the following
conditions:
\begin{itemize}
\item[(i)]
for every space $X$ and each discrete space $D$, the space $X$ has property
$\mathscr{E}$ if and only if $X\times D$ has property $\mathscr{E}$;
\item[(ii)]
if $X\in\mathscr{C}$ and $D$ is a discrete space, then $X\times D\in\mathscr{C}$.
\end{itemize}
Assume that
$G^\kappa$-equivalence preserves property $\mathscr{E}$ within the class $\mathscr{C}$ for
some cardinal $\kappa\ge 1$. Then $G^\tau$-equivalence preserves property $\mathscr{E}$ within
the class $\mathscr{C}$ for each cardinal $\tau\ge 1$.
\end{proposition}
\begin{proof}
The assumption of our
proposition
yields $\mathbf{I}(G,\mathscr{E},\mathscr{C})\not=\emptyset$, and so
$1\in \mathbf{I}(G,\mathscr{E},\mathscr{C})$ by Lemma
\ref{non:empty:implies:1}. Applying the assumption of our
proposition
 and Lemma \ref{1:implies:tau} to each cardinal $\tau\ge 1$, we conclude that
$\mathbf{I}(G,\mathscr{E},\mathscr{C})=\{\tau\ge 1: \tau$ is a
cardinal$\}$.
\end{proof}

The above proposition is applicable to many ``local'' properties.
\begin{corollary}
Let $\mathscr{C}$ be a class of spaces such that $X\times
D\in\mathscr{C}$ whenever $X\in\mathscr{C}$ and $D$ is a discrete
space. Let $\mathscr{E}$ be one of the following properties:
metrizability, paracompactness, weak paracompactness, local
compactness, first countability, countable tightness,
Frech\'et-Urysohn property, sequentiality,
total disconnectedness,
(property of having a given) covering dimension $\dim$, large inductive dimension $\mathrm{Ind}$,
small inductive dimension $\mathrm{ind}$.
Assume also that $G$ is a
topological group
such that $G^\kappa$-equivalence preserves property $\mathscr{E}$
within the class $\mathscr{C}$
for some cardinal $\kappa\ge 1$.
Then for every cardinal
$\tau\ge1$,
 $G^\tau$-equivalence preserves
property $\mathscr{E}$ within the class $\mathscr{C}$.
\end{corollary}

The list of properties $\mathscr{E}$ in the above corollary can be
easily extended.

\begin{proposition}
\label{findissum}
Let $G$ be a topological group and $\kappa$ an infinite cardinal.
Assume that $\mathscr{E}$ is
a topological property and $\mathscr{C}$ is a class of spaces such that:
\begin{itemize}
\item[(i)]
for every space $X$, the
space
 $X\times D_\tau$
has property $\mathscr{E}$ if and only if $X$ has property $\mathscr{E}$ and
$1\le\tau<\kappa$,
\item[(ii)] if $X\in\mathscr{C}$ and $D$ is a discrete space, then $X\times D\in\mathscr{C}$, and
\item[(iii)]
there exists at least one space $X\in\mathscr{C}$ satisfying property $\mathscr{E}$.
\end{itemize}
Then either $\mathbf{I}(G,\mathscr{E},\mathscr{C})=\emptyset$, or
\begin{equation}
\label{eq:G:E}
\mathbf{I}(G,\mathscr{E},\mathscr{C})=\{\tau: \tau
\mbox{ is a cardinal satisfying }
1\le \tau<\kappa\}.
\end{equation}
\end{proposition}

\begin{proof}
Let
$\tau\ge\kappa$.
In particular, $\tau$ is infinite.
By item (iii), there exists a space $X\in\mathscr{C}$ satisfying property $\mathscr{E}$.
Then $X\times D_\tau\in\mathscr{C}$ by item (ii).
Since $\tau\ge\kappa$, from item (i) we conclude that  $X\times D_\tau$ does
not have property $\mathscr{E}$. Since $X$ and $X\times D_\tau$ are
$G^\tau$-equivalent by Proposition \ref{group.of.countable.power}(ii), it follows
that $\tau\not\in \mathbf{I}(G,\mathscr{E},\mathscr{C})$.
We have proved that
\begin{equation}
\label{eq:G:E:inclusion}
\mathbf{I}(G,\mathscr{E},\mathscr{C})\subseteq\{\tau: \tau
\mbox{ is a cardinal satisfying }
1\le \tau<\kappa\}.
\end{equation}

Assume now that $\mathbf{I}(G,\mathscr{E},\mathscr{C})\not=\emptyset$.
Then $1\in \mathbf{I}(G,\mathscr{E},\mathscr{C})$
by Lemma \ref{non:empty:implies:1}.
If
$\tau$ is a cardinal satisfying
 $1\le \tau<\kappa$, then from (i),
(ii) and Lemma \ref{1:implies:tau} we get  $\tau\in
\mathbf{I}(G,\mathscr{E},\mathscr{C})$. Together with
\eqref{eq:G:E:inclusion} this proves \eqref{eq:G:E}.
\end{proof}

\begin{corollary}
\label{preservation:by:powers} Let $G$, $\kappa$, $\mathscr{E}$ and
$\mathscr{C}$ be as in the assumption of Proposition \ref{findissum}.
Assume also that
 $G^\sigma$-equivalence
preserves the property $\mathscr{E}$
within the class $\mathscr{C}$
for some cardinal $\sigma\ge 1$.
Then the following conditions are
equivalent:
\begin{itemize}
\item[(i)]
$G^\tau$-equivalence preserves the property $\mathscr{E}$ within the class $\mathscr{C}$;
\item[(ii)]
$1\le \tau<\kappa$.
\end{itemize}
\end{corollary}
\begin{proof}
The assumption of our corollary yields
$\sigma\in\mathbf{I}(G,\mathscr{E},\mathscr{C})\not=\emptyset$, so
\eqref{eq:G:E} holds by Proposition \ref{findissum}. It remains only
to show that $G^0$-equivalence does not preserve the property
$\mathscr{E}$ within the class $\mathscr{C}$. Since $G^0=\{e\}$ is
the trivial group, so is $\Cp{Y}{G^0}$ for every space $Y$. Thus,
 any two spaces are $G^0$-equivalent. In particular, $X\eq{G^0}X\times
 D_\kappa$, where $X\in\mathscr{C}$ is
the space from item (iii) of Proposition \ref{findissum}. Then $X\times
D_\kappa\in\mathscr{C}$  by Proposition \ref{findissum}(ii), while
$X\times D_\kappa$ does not have property  $\mathscr{E}$ by
Proposition \ref{findissum}(i). We conclude that $G^0$-equivalence
does not preserve the property $\mathscr{E}$ within the class
$\mathscr{C}$.
\end{proof}

\begin{corollary}
\label{productOK} Let $G$ be a topological group, $m\ge 1$
an
integer number
and
$\mathscr{C}$
a class of spaces satisfying the condition from item (ii) of
Proposition \ref{findissum}.
If $G^m$-equivalence preserves compactness (countable compactness,
pseudocompactness, $\sigma$-compactness)
within the class  $\mathscr{C}$,
then $G^k$-equivalence preserves the corresponding property
 within the class  $\mathscr{C}$
for every $k\in\mathbb{N}\setminus\{0\}$.
\end{corollary}
\begin{proof}
It suffices to note that $\kappa=\omega$, $\mathscr{E}$ and $\mathscr{C}$ satisfy the assumptions of Proposition \ref{findissum},
where $\mathscr{E}$ is one of the four properties listed in the statement of our corollary. Since
$m\in \mathbf{I}(G,\mathscr{E},\mathscr{C})\not=\emptyset$,
the conclusion of our corollary follows from that of
Proposition \ref{findissum}.
\end{proof}

\section{When does $(G\times H)$-equivalence imply both $G$-equivalence and $H$-equivalence?}
\label{section:3}

Let $G$ and $H$ be topological groups. If two spaces are both
$G$-equivalent and $H$-equivalent, then they are also $(G\times
H)$-equivalent (Proposition \ref{product:equivalence}).
In Example \ref{GxH:does:not:imply:G} we exhibit topological groups
$G$ and $H$ such that $(G\times H)$-equivalence implies neither $G$-equivalence, nor
$H$-equivalence,
thereby demonstrating that the converse implication
fails in general.
 In this
section we will prove that this implication holds {\em for
sufficiently different\/} topological groups $G$ and $H$; see
Theorems \ref{thm:4.4} and \ref{thm:4.5}.
These
results
turn out to be useful
in constructing
numerous
examples in Sections \ref{seccomp} and \ref{secdis}.

\begin{lemma}
\label{splitting:lemma}
Suppose that $X$ and $Y$ are spaces, $G$ and $H$ are topological groups and
\begin{equation}
\label{isom:eq}
\varphi:\Cp{X}{G}\times\Cp{X}{H}\rightarrow\Cp{Y}{G}\times\Cp{Y}{H}
\mbox{ is a topological isomorphism }
\end{equation}
satisfying
\begin{equation}
\label{G:1:H:1}
\varphi\left(\Cp{X}{G}\times \{1_{X}\}\right)
=
\Cp{Y}{G}\times \{1_{Y}\},
\end{equation}
where $1_{X}$ and $1_{Y}$ denote the identity elements of $\Cp{X}{H}$ and $\Cp{Y}{H}$, respectively.
Then $X$ and $Y$ are both $G$-equivalent and $H$-equivalent.
\end{lemma}
\begin{proof}
From \eqref{isom:eq} and
\eqref{G:1:H:1}
we get
$\Cp{X}{G}\cong \Cp{Y}{G}$. This proves that $X$ and $Y$ are $G$-equivalent.
Applying \eqref{isom:eq} and
\eqref{G:1:H:1} once again, one easily obtains that
\begin{align*}
\Cp{X}{H}
&\cong
\left(\Cp{X}{G}\times\Cp{X}{H}\right)
 / \left(\Cp{X}{G}\times \{1_{X}\}\right)
\\
&
\cong
\left(\Cp{Y}{G}\times\Cp{Y}{H}\right)
 / \left(\Cp{Y}{G}\times \{1_{Y}\}\right) \cong \Cp{Y}{H},
\end{align*}
which proves that $X$ and $Y$ are $H$-equivalent.
\end{proof}

\begin{theorem}
\label{thm:4.4} Let $G$ and $H$ be topological groups satisfying one
of the following two conditions:
\begin{itemize}
\item[(i)]
$G$ is pathwise connected and $H$ is hereditarily
 disconnected;
\item[(ii)]
 $G$ is a precompact group and $H$ is a topological group without
 nontrivial precompact subgroups.
\end{itemize}
If two spaces are $(G\times H)$-equivalent, then they are both
$G$-equivalent and $H$-equivalent.
\end{theorem}
\begin{proof}
Assume that $X$ and $Y$ are $(G\times H)$-equivalent spaces. Then
$$
\Cp{X}{G}\times\Cp{X}{H}\cong \Cp{X}{G\times H}\cong \Cp{Y}{G\times
H} \cong \Cp{Y}{G}\times\Cp{Y}{H},
$$
so we can fix $\varphi$ satisfying \eqref{isom:eq}. We continue
using notation from Lemma \ref{splitting:lemma}. There are two
cases.

{\sl Case 1\/}. {\it Item (i) holds\/}. Suppose that $n\in\N$,
$g_1,\ldots g_n\in G$ and $x_1,\ldots,x_n\in X$ are pairwise
distinct. Fix $i=1,\dots,n$. Since $G$ is pathwise connected, there
exists a continuous map $\varphi_i: [0,1]\to G$ such that
\begin{equation}
\label{eq:varphi:i} \varphi_i(0)=e \mbox{  and } \varphi_i(1)=g_i.
\end{equation}
Let $\psi_i: X\to [0,1]$ be a continuous function such that
\begin{equation}
\label{eq:psi:i} \psi_i(x_i)=1 \mbox{ and } \psi_i(x_j)=0 \mbox{ for
every } j\in\{1,\dots,n\} \mbox{ with } j\not=i.
\end{equation}

Let $\varphi:[0,1]\to \Cp{X}{G}$ be the map which assigns to every
$t\in [0,1]$ the function $\varphi(t)\in\Cp{X}{G}$ defined by
\begin{equation}
\label{eq:final:oath}
\varphi(t)(x)=\prod_{i=1}^n\varphi_i(t\psi_i(x)) \mbox{ for } x\in
X.
\end{equation}
From \eqref{eq:varphi:i} and \eqref{eq:final:oath} we conclude that
$\varphi(0)$ is the identity element of $\Cp{X}{G}$. From
\eqref{eq:varphi:i}, \eqref{eq:psi:i} and \eqref{eq:final:oath} we
conclude that $h(x_i)=g_i$ for every integer $i$ with $1\le i\le n$,
where $h=\varphi(1)$. One can easily check that $\varphi$ is
continuous. This argument proves that $\Cp{X}{G}$ is connected.

Since $H$ is hereditarily disconnected, so is $H^X$. Since
$\Cp{X}{H}\subseteq H^X$, $\Cp{X}{H}$ is hereditarily disconnected
as well. It follows that
$c(\Cp{X}{G}\times\Cp{X}{H})=\Cp{X}{G}\times\{1_X\}$. Similarly,
$c(\Cp{Y}{G}\times\Cp{Y}{H})=\Cp{Y}{G}\times \{1_Y\}$. Since
$\varphi$ is a topological isomorphism, we obtain \eqref{G:1:H:1}.

{\sl Case 2\/}. {\it Item (ii) holds\/}. For $y\in Y$ let $\varpi_y:
\Cp{Y}{G}\times\Cp{Y}{H}\to H$ be the continuous homomorphism
defined by $\varpi_y(g,h)=h(y)$ for $(g,h)\in
\Cp{Y}{G}\times\Cp{Y}{H}$.

Since $G$ is precompact, so is $G^X$. Being a subgroup of the
precompact group $G^X$, the group $\Cp{X}{G}$ is precompact as well.
Being an image of the precompact group under a continuous group
homomorphism, $H_y=\varpi_y(\varphi\left(\Cp{X}{G}\times
\{1_X\}\right))$ is a precompact subgroup of $H$ for every $y\in Y$.
By our assumption, each $H_y$ must be the trivial subgroup of $H$,
which yields the inclusion $ \varphi\left(\Cp{X}{G}\times
\{1_X\}\right) \subseteq \Cp{Y}{G}\times \{1_Y\} $. By the symmetry, the
inclusion $ \varphi^{-1}\left(\Cp{Y}{G}\times \{1_Y\}\right)
\subseteq \Cp{X}{G}\times \{1_X\} $ holds as well.
This proves \eqref{G:1:H:1}.

In both cases the conclusion follows from Lemma
\ref{splitting:lemma}.
\end{proof}

Let $G$ be a group.
Recall that $g\in G$ is called a {\em torsion element of $G$} if there exists some
$n\in\mathbb{N}\setminus\{0\}$ such that $g^n=e$. The
subset
of all torsion
elements of $G$ is called the {\em torsion part of $G$\/} and denoted
by \tor{G}. If \tor{G}=$\{e\}$, then $G$ is called {\em torsion-free}.
For a given $n\in\mathbb{N}$, let
$G^{(n)}=\{g^n:g\in G\}$.

\begin{theorem}
\label{thm:4.5}
Let $G$ and $H$ be topological groups satisfying one of the following two conditions:
\begin{itemize}
\item[(i)]
\tor{{G}} is dense in
$G$
and $\hat{H}$ is torsion-free;
\item[(ii)]
there exists $n\in\mathbb{N}$ such that
$\hat{G}^{(n)}=\hat{G}$ and
${H}^{(n)}=\{e\}$.
\end{itemize}
If $G^{\star\star}$-regular spaces are $(G\times H)$-equivalent,
then they are both $G$-equivalent and $H$-equivalent.
\end{theorem}
\begin{proof}
Assume that $G^{\star\star}$-regular spaces $X$ and $Y$ are $(G\times H)$-equivalent.
Arguing as in the beginning of the proof of Theorem \ref{thm:4.4}, we can fix $\varphi$ satisfying
\eqref{isom:eq}.
For typographical reasons, define
$$
G(X)=\widehat{\Cp{X}{G}},\ \
G(Y)=\widehat{\Cp{Y}{G}},\ \
H(X)=\widehat{\Cp{X}{H}}
\ \mbox{ and }\
H(Y)=\widehat{\Cp{Y}{H}}.
$$
Let $ \Phi:G(X)\times H(X)\rightarrow G(Y)\times H(Y) $ be the
(unique) topological isomorphism extending $\varphi$. We claim that
\begin{equation}
\label{eq:Phi}
\Phi\left(G(X)\times \{e_{H(X)}\}\right)
=
G(Y)\times \{e_{H(Y)}\}.
\end{equation}
Since $X$
is
$G^{\star\star}$-regular,
one can easily see that
$\Cp{X}{G}$ is dense in $G^X$.
Therefore,
\begin{equation}
\label{computation:of:G(X)}
G(X)=\widehat{\Cp{X}{G}}=\widehat{G^X}=\widehat{G}^X.
\end{equation}
We need to consider two cases.

{\sl Case 1\/}. {\it Item (i) holds\/}. Since \tor{G} is dense in
$G$, and the latter group is dense in $\hat{G}$, we conclude that
$tor(\hat{G}^X)$ is dense in $\hat{G}^X$. Combining this with
\eqref{computation:of:G(X)}, we conclude that $tor(G(X))$ is dense
in $G(X)$. Since $\Phi$ is a topological isomorphism, it follows
that
\begin{equation}
\label{eq:4.9}
\Phi\left(tor(G(X))\times \{e_{H(X)}\}\right)
\mbox{ is dense in }
\Phi\left(G(X)\times \{e_{H(X)}\}\right)
\end{equation}
and
\begin{equation}
\label{eq:4.10}
\Phi\left(tor(G(X))\times \{e_{H(X)}\}\right)
\subseteq
tor\left(G(Y)\times H(Y)\right).
\end{equation}
Since
$H(Y)=\widehat{\Cp{Y}{H}}\subseteq
\widehat{H}^Y$
and $\widehat{H}$ is torsion-free,
$H(Y)$ must be torsion-free as well.
In particular,
\begin{equation}
\label{eq:4.11}
tor\left(G(Y)\times H(Y)\right)
\subseteq G(Y)\times \{e_{H(Y)}\}.
\end{equation}
Since the latter set is closed in $G(Y)\times H(Y)$, from
\eqref{eq:4.9}, \eqref{eq:4.10} and \eqref{eq:4.11} one concludes
that
\begin{equation}
\label{eq:4.12}
\Phi\left(G(X)\times \{e_{H(X)}\}\right)
\subseteq
G(Y)\times \{e_{H(Y)}\}.
\end{equation}
Applying the same arguments to the inverse map
$\Phi^{-1}$ of $\Phi$,
we get
\begin{equation}
\label{eq:4.13}
\Phi^{-1}\left(G(Y)\times \{e_{H(Y)}\}\right)
\subseteq
G(X)\times \{e_{H(X)}\}.
\end{equation}

{\sl Case 2\/}. {\it Item (ii) holds\/}. Fix $n\in\mathbb{N}$ as in
item (ii). Choose $g\in G(X)$ arbitrarily. From
$\hat{G}^{(n)}=\hat{G}$ one gets
$\left(\hat{G}^X\right)^{(n)}=\hat{G}^X$. Combining this with
\eqref{computation:of:G(X)}, we obtain that $G(X)^{(n)}=G(X)$. Hence
there exists $g_0\in G(X)$ such that $g_0^n=g$. Let
$\Phi\left(g_0,e_{H(X)}\right)=(f,h)$.

Since $H^{(n)}=\{e_H\}$, by the ``principle of extending of equations'',
$\hat{H}^{(n)}=\{e_{\hat{H}}\}$ holds as well.
Since $h\in H(Y)\subseteq \hat{H}^Y$,
we conclude that $h^n=e_{H(Y)}$.
Thus,
$$
\Phi\left(g,e_{H(X)}\right)=
\Phi\left(\left(g_0,e_{H(X)}\right)^n\right)
=
\Phi\left(g_0,e_{H(X)}\right)^n
=
(f,h)^n=(f^n,h^n)
=
\left(f^n,e_{H(Y)}\right).
$$
This proves the inclusion \eqref{eq:4.12}. Applying the same
arguments to the inverse map $\Phi^{-1}$ of $\Phi$, we get the
inclusion \eqref{eq:4.13}.

Going back to the common proof, note that \eqref{eq:4.12} and
\eqref{eq:4.13} yield \eqref{eq:Phi}.
Furthermore, since $\Phi$ extends the isomorphism
$\varphi$, from \eqref{eq:Phi} one gets \eqref{G:1:H:1}, and now the
application of Lemma \ref{splitting:lemma} finishes the proof.
\end{proof}

Recall that, for a prime number $p$, a group $G$ {\em has exponent $p$\/} if
$G^{(p)}=\{e_G\}$ and $G\not=\{e_G\}$.

\begin{corollary}
Let $p$ and $q$ be distinct prime numbers.
Suppose that a topological group $G$ has exponent $p$
and a topological group $H$ has exponent $q$.
Suppose also that spaces $X$ and $Y$ are $G^{\star\star}$-regular
and $(G\times H)$-equivalent.
Then $X$ and $Y$ are both $G$-equivalent and $H$-equivalent.
\end{corollary}
\begin{proof}
Since $G$ has exponent $p$,
the ``principle of extending of equations'' implies
$\hat{G}^{(p)}=\{e_{\hat{G}}\}$. That is, $\hat{G}$ has exponent $p$ as
well.
One can easily see that this yields $\hat{G}^{(q)}=\hat{G}$.
Since
${H}^{(q)}=\{e_H\}$, the conclusion of our corollary follows from
Theorem \ref{thm:4.5}(ii).
\end{proof}

\begin{example}
\rm{Let $(\mathbb{C}\setminus\{0\},\cdot)$ denote the multiplicative
group of all nonzero complex numbers with its usual topology. Then
spaces $X$ and $Y$ are $(\mathbb{C}\setminus\{0\},\cdot)$-equivalent
if and only if they are both $\mathbb{R}$-equivalent and
$\mathbb{T}$-equivalent. Indeed, it suffices to realize that
$\mathbb{C}\setminus\{0\}\cong\mathbb{T}\times\mathbb{R}$. The rest
follows from
Proposition \ref{product:equivalence},
combined with either
Theorem \ref{thm:4.4}(ii), or Theorem \ref{thm:4.5}(i).}
\end{example}

\section{\IANS\  and NSS groups}
\label{section:4}

\begin{definition}
{\rm
 We say that a subset $A$
of a topological group $G$ is {\em \ap\ in
$G$\/} provided that, for every injection $a:\mathbb{N}\rightarrow A$ and each
mapping $z:\mathbb{N}\rightarrow\mathbb{Z}$,
the sequence
\begin{equation}
\label{defining:sequence:for:infinite:products}
\left\{\prod_{n=0}^{k}a(n)^{z(n)}:k\in\N\right\}
\end{equation}
 of elements of
$G$ converges to some $g\in G$. In such a case we will also say that the
(infinite) {\em product $\prod_{n=0}^\infty a(n)^{z(n)}$ converges to $g$\/}
and write
\begin{equation}
\label{g:is:infinite:product}
g=\prod_{n=0}^\infty a(n)^{z(n)}.
\end{equation}
 }
\end{definition}

The proofs of the next three lemmas are straightforward.
\begin{lemma}
\label{subset:of:productive:is:prodictive}
\label{a.p.ispreserved}
\begin{itemize}
\item[(i)]
A subset of an \ap\ set is \ap.
\item[(ii)]
Let $\phi:G\rightarrow H$ be a continuous homomorphism between
topological groups $G$ and $H$. If a set $A\subseteq G$ is  \ap\  in $G$, then
$\phi(A)$ is \ap\  in $H$.
\end{itemize}
\end{lemma}

\begin{lemma}
\label{basic:properties:of:ap}
Let $H$ be a subgroup of a topological group $G$ and $A\subseteq H$.
\begin{itemize}
\item[(i)]
If $A$ is \ap\ in $H$, then it is \ap\ in $G$ as well.
\item[(ii)]
If $H$ is sequentially closed in $G$, then $A$ is \ap\ in $H$ if and only if $A$ is \ap\ in $G$.
\end{itemize}
\end{lemma}

Item (i) of the next lemma
gives a typical example of an infinite
\ap\ set, while item (ii) shows that ``sequentially closed'' cannot be omitted in Lemma \ref{basic:properties:of:ap}(ii).
\begin{lemma}
\label{algebraically:productive:subsets:in:products} Let $\{G_i\in
I\}$ be
an infinite
family consisting of nontrivial topological groups $G_i$.
Let $G=\prod_{i\in I} G_i$
and
$$
H=\{g\in G: \mbox{ the set } \{i\in I: g(i)\not=e_{G_i}\} \mbox{ is finite}\}.
$$
For each $i\in I$ choose $g_i\in
G\setminus\{e\}$ such that $g_i(j)= e_{G_j}$ for every $j\in
I\setminus\{i\}$, and consider the infinite set $A=\{g_i:i\in I\}\subseteq H$.
Then:
\begin{itemize}
\item[(i)]
$A$
 is
 \ap\
in
$G$, but
\item[(ii)]
$A$ is not
\ap\
in
$H$.
\end{itemize}
\end{lemma}

\begin{definition}
{\rm We say that a topological group $G$ is {\em $\IANS$\/} (an abbreviation for ``Trivially
Absolutely Productive'') if every \ap\ set in
$G$ is finite. }
\end{definition}

\begin{proposition}\label{h}
\begin{enumerate}
\item[(i)]The class of all \IANS\ groups is closed under taking finite
products and subgroups.
\item[(ii)] Let $G_i$ be a nontrivial topological group for
every $i\in\mathbb{N}$. Then $G=\prod_{i=0}^\infty G_i$ is not \IANS.
\item[(iii)]
A \IANS\ group does not contain any subgroup topologically isomorphic to a Cartesian
product
$\prod_{i=0}^\infty G_i$
of
nontrivial topological
groups $G_i$.
\end{enumerate}
\end{proposition}
\begin{proof}
Item (i) follows from Lemmas \ref{a.p.ispreserved}(ii) and
\ref{basic:properties:of:ap}(i),
 item (ii) follows from Lemma
\ref{algebraically:productive:subsets:in:products}(i), and item (iii)
follows from items (i) and (ii).
\end{proof}

\begin{remark}
{\rm The converse of Proposition \ref{h}(iii) does not hold in
general. Indeed, the group $\mathbb{Z}_p$ of $p$-adic integers does
not contain any subgroup topologically isomorphic to an infinite
product of nontrivial groups, yet $\mathbb{Z}_p$ is not \IANS\
\cite{DSS}. }
\end{remark}

Recall that a topological group $G$ is an {\em NSS group\/}, or has
an {\em NSS property\/} (an abbreviation for ``no small subgroups'')
if $G$ has an open
neighborhood of
the
identity containing no
nontrivial subgroups of $G$. The following lemma provides a simple
reformulation of the NSS property.
\begin{lemma}\label{NSS}
\label{charNSS} Let $G$ be a topological group. Then the following conditions
are equivalent:
\begin{enumerate}
 \item[(i)] $G$ is an NSS group;
 \item[(ii)] there exists an open neighborhood $U$ of the identity $e$ of $G$ such that for
every $g\in G\setminus\{e\}$ and each $f,h\in G$
one can find
$z\in\mathbb{Z}$ with  $hg^z\not\in fU$.
\end{enumerate}
\end{lemma}
\begin{proof}
(i)$\Rightarrow$(ii) Let $V$ be a neighborhood of the identity $e$
witnessing that $G$ is NSS. Choose a
neighborhood $U$ of $e$
with $U^{-1}U\subseteq V$.
Let $g\in
G\setminus\{e\}$ and $f,h\in G$ be arbitrary.
If $f^{-1}h\not\in U$,
then $f^{-1}hg^0=f^{-1}h\not\in U$, and consequently $hg^0\not\in
fU$, so $z=0$ works.
Suppose now that
$f^{-1}h\in U$. Then $h^{-1}fU\subseteq U^{-1}U\subseteq V$.
By the choice of $V$, we can find
$z\in\mathbb{Z}$ such that $g^z\not\in V$ (otherwise $V$
would contain the nontrivial cyclic subgroup generated be $g$).
In particular,
$g^z\not\in
h^{-1}fU$, and so $hg^z\not\in fU$.

(ii)$\Rightarrow$(i)
Applying (ii) with $f=h=e$, we conclude that for every $g\in
G\setminus\{e\}$ there exists $z\in\mathbb{Z}$ such that $g^z\not\in
U$. This means that $U$ is an open neighborhood of
$e$
which
contains no nontrivial subgroup. Thus $G$ is NSS.
\end{proof}

\begin{theorem}\label{NSSjeTAP}
An NSS group is \IANS.
\end{theorem}
\begin{proof}
Assume that $A$ is an infinite subset of an NSS group $G$. We must
show that $A$ is not \ap\ in $G$.
Fix an injection $a:\mathbb{N}\to A\setminus\{e\}$.
The
subgroup $H$ of $G$ generated by $a(\mathbb{N})$ is countable,
so we can choose a
map $f:\mathbb{N}\to H$ such that
$f^{-1}(h)$ is infinite for every $h\in H$.
Since $G$ is NSS, we can fix $U$ satisfying item (ii) of Lemma
\ref{NSS}.
We are going to define a map $z:\mathbb{N}\to \mathbb{Z}$ such that
\begin{equation*}
\label{eq:5.10:16}
\prod_{i=0}^{k}a(i)^{z(i)}\not\in f(k) U
\tag{*$_k$}
\end{equation*}
holds for every $k\in\mathbb{N}$. Since $a(0)\not=e$, applying item
(ii) of Lemma \ref{NSS} to $g=a(0)$, $f=f(0)$ and $h=e$,
we can choose $z(0)\in\mathbb{Z}$ satisfying (*$_0$).
For
$n\in\mathbb{N}\setminus\{0\}$, assume that $z(j)\in\mathbb{Z}$
satisfying (*$_j$) has already been defined for all $j<n$. Since
$a(n)\not=e$, applying item (ii) of Lemma \ref{NSS} to $f=f(n)$,
$g=a(n)$ and $h=\prod_{i=0}^{n-1}a(i)^{z(i)}$ we can choose
$z(n)\in\mathbb{Z}$
 satisfying (*$_n$). This finishes the inductive
construction.

Suppose now that \eqref{g:is:infinite:product}
holds
for some $g\in G$.
Since $\prod_{i=0}^{k}a(i)^{z(i)}\in H$ for every $k\in\mathbb{N}$,
$g$ must belong to the sequential closure of $H$.
In particular, $g\in\overline{H}\subseteq HU$, and so
$g\in hU$ for some $h\in H$.
On the other hand, \eqref{eq:5.10:16} holds for every $k\in\mathbb{N}$, which gives
$$
f^{-1}(h)\subseteq \left\{k\in\mathbb{N}: \prod_{i=0}^{k}a(i)^{z(i)}
\not\in hU\right\}.
$$
Since the set $f^{-1}(h)$ is infinite and $g\in hU$, we conclude that the sequence
\eqref{defining:sequence:for:infinite:products}
cannot converge to $g$, in contradiction with
\eqref{g:is:infinite:product}. This proves that $A$ is not \ap\ in $G$.
\end{proof}

\begin{proposition}
\label{psc:TAP:not:NSS}
\begin{itemize}
\item[(i)]
Every topological group
without nontrivial convergent sequences is $\IANS$.
\item[(ii)]
An infinite pseudocompact group
without
nontrivial convergent sequences is \IANS\ but not NSS.
\end{itemize}
\end{proposition}
\begin{proof}
(i)
Let $A$ be an infinite subset of a topological group $G$
without
nontrivial convergent sequences.
Fix an injection $a:\mathbb{N}\to A$. By induction on
$m\in\mathbb{N}$ one can easily choose an increasing sequence
$\{n_m:m\in\mathbb{N}\}\subseteq \mathbb{N}$
 such that
\begin{equation}
\label{eq:5.8:16}
\prod_{i=0}^m a(n_i)\not\in
\left\{
\prod_{i=0}^k a(n_i): k\in\mathbb{N}, k<m
\right\}.
\end{equation}
Define $z:\mathbb{N}\to\mathbb{Z}$ by $z(n)=1$ if $n\in
\{n_m:m\in\mathbb{N}\}$ and $z_n=0$ otherwise.  It follows from
\eqref{eq:5.8:16} that the sequence
\eqref{defining:sequence:for:infinite:products} is nontrivial, and
so it cannot converge by our assumption. Thus, $A$ is not \ap\ in $G$.

(ii) Assume, in addition, that $G$ is pseudocompact and infinite.
Suppose that $G$ is NSS, and choose an open neighborhood $U$ of the
identity $e$ that contains no nontrivial subgroups of $G$. Starting with
$U_0=U$, choose a sequence
$\{U_n: n\in\mathbb{N}\}$
 of open neighborhoods of $e$ such that
$U_{n+1}^{-1} U_{n+1}\subseteq U_n$ for every $n\in\mathbb{N}$. Then
$H=\bigcap_{n=0}^\infty U_n\subseteq U_0=U$ is a subgroup of $G$,
which gives $H=\{e\}$. Since $G$ is pseudocompact, we conclude that
$G$ must be metrizable and hence compact. Being infinite, $G$ must
contain a nontrivial convergent sequence, a contradiction.
\end{proof}

\begin{remark}
\label{compact-like:TAP:vs:NSS}
\rm{
\begin{itemize}
\item[(i)] There is an infinite pseudocompact Abelian
group without nontrivial convergent sequences \cite{Sirota}.
Therefore, from Proposition \ref{psc:TAP:not:NSS} we
conclude that {\em there exists a pseudocompact Abelian \IANS\ group
that is not NSS\/}.
\item[(ii)] There are consistent examples of infinite countably compact
Abelian groups without nontrivial convergent sequences; see
\cite{DS} for references. Applying Proposition
\ref{psc:TAP:not:NSS}, we conclude that {\em the existence of a
countably compact Abelian \IANS\ group that is not NSS is consistent
with ZFC\/}.
\item[(iii)]
Theorem \ref{NSSjeTAP} can sometimes be reversed.
Indeed, it has been proved recently in
\cite{DSS} that {\em a
locally compact
 \IANS\ group $G$ is NSS\/}.
Moreover, {\em a totally disconnected compact \IANS\ group is finite\/} \cite{DSS}.
\item[(iv)] It is proved in \cite{DSS} that {\em a $\sigma$-compact complete Abelian \IANS\ group need not be NSS\/}.
\end{itemize}
}
\end{remark}

We refer the reader to Proposition \ref{Cp:TAP:not:NSS}(ii) for other examples of \IANS\
groups that are not NSS.

\begin{remark}
{\rm A short alternative proof of Theorem \ref{NSSjeTAP} has been given recently in \cite{DT}.}
\end{remark}

\section{A group-theoretic proof that $l$-equivalence preserves pseudocompactness}
\label{section:5}

\begin{lemma}
\label{konecne} Assume that $X$ is a space, $G$ is
a \IANS\
group
and
$A$ is a subset of $\Cp{X}{G}$.
If $A$ is \ap\ in $\Cp{X}{G}$,
then the set $\{f\in A: f(x)\not=e\}$ is finite  for every $x\in X$.
\end{lemma}
\begin{proof}
Suppose, by the way of contradiction, that there exists $x\in X$
such that $f(x)\neq e$  for infinitely many $f\in A$. For each
$f\in\Cp{X}{G}$ define $\pi_x(f)=f(x)$. Since $\pi_x:\Cp{X}{G}\to G$
is a continuous homomorphism,
$\pi_x(A)$ is \ap\ in $G$
by Lemma \ref{a.p.ispreserved}(ii).
Since $G$ is \IANS, $\pi_x(A)$ must be finite, so
by the pigeon hole principle, there exists $g\in
G\setminus\{e\}$ and an infinite set $\{f_{i}:i\in\mathbb{N}\}\subseteq A$
such that $f_i(x)=g$ for every $i\in\mathbb{N}$.
Since
the
sequence $\left\{\prod_{i=0}^k f_i^{(-1)^i}(x):k\in\mathbb{N}\right\}$
alternates between $e$ and $g\not=e$,
the product $\prod_{i=0}^\infty f_i^{(-1)^i}$ does not exists.
This contradicts the fact that $A$ is \ap.
\end{proof}

\begin{lemma}\label{selecting:points}
Let $X$ be a space, $G$ a \IANS\ group and $A$ an infinite \ap\ subset of $\Cp{X}{G}$.
Then there exist two (necessarily faithfully indexed)
sequences $\{x_i:i\in\mathbb{N}\}\subseteq X$
and
$\{f_i:i\in\mathbb{N}\}\subseteq A$
such that
\begin{equation}\label{eq2}
 \mbox{$f_{i}(x_{i})\neq e$ and $f_{i}(x_{j})=e$ whenever
$i,j\in\mathbb{N}$ and $j<i$.}
\end{equation}
\end{lemma}
\begin{proof}
We use induction on $i\in\mathbb{N}$. First,
choose $f_0\in A$ and $x_0\in X$ such that $f_0(x_0)\neq e$.
Let $n\in\mathbb{N}\setminus\{0\}$, and suppose that
$\{x_0,x_1,\dots,x_{n-1}\}\subseteq X$ and
$\{f_0,f_1,\dots,f_{n-1}\}\subseteq A$
have
 already been
selected so that  $f_i(x_i)\not=e$ and $f_i(x_j)=e$ whenever $i,j
\in\mathbb{N}$ and $j <i \le n-1$.
  The set
$B_n=\bigcup_{i=0}^{n-1}\{f\in A:f(x_i)\neq e\}$ is finite by Lemma
\ref{konecne}, and hence there exists $f_n\in A\setminus
B_n\not=\emptyset$. Without loss of generality, $f_n$ is not the
identity element of $\Cp{X}{G}$, and so $f_n(x_n)\not=e$ for some
$x_n\in X$.
\end{proof}

\begin{theorem} \label{Rpseudo}
A space $X$ is pseudocompact if and only if
$\Cp{X}{\mathbb{R}}$ has
the
\IANS\ property.
\end{theorem}
\begin{proof}
Since the group $\Cp{X}{\mathbb{R}}$ is Abelian, in this proof we shall
use the additive notation.

To prove the ``if'' part, suppose that $X$ is not pseudocompact. Fix
an infinite discrete family $\mathscr{U}=\{U_i:i\in\mathbb{N}\}$ of
non-empty open subsets of $X$. For each $i\in\mathbb{N}$ choose
$x_i\in U_i$ and $f_i\in \Cp{X}{\mathbb{R}}$ such that $f(x_i)\not =
0$ and $f(X\setminus U_i)\subseteq \{0\}$. Clearly,
$A=\{f_i:i\in\mathbb{N}\}$ is faithfully indexed (and thus
infinite). Let $s:\mathbb{N}\to \mathbb{N}$ be an injection and
$z:\mathbb{N}\to\mathbb{Z}$ a map. Since $\mathscr{U}$ is discrete,
$f=\sum_{i=0}^{\infty} z(i) f_{s(i)}\in \Cp{X}{\mathbb{R}}$
and
$f=\lim_{k\to\infty} \sum_{i=0}^{k} z(i) f_{s(i)}$. This
shows that $A$ is \ap\ in $\Cp{X}{\mathbb{R}}$, and so
$\Cp{X}{\mathbb{R}}$ is not \IANS.

Being an NSS group, $\mathbb{R}$ has
the
\IANS\ property by Theorem
 \ref{NSSjeTAP}. A simpler direct proof can be obtained as follows.
Let $A$ be an infinite subset of $\mathbb{R}$.
Fix an injection
$a: \mathbb{N}\to A\setminus\{0\}$. By induction on $k\in\mathbb{N}$
choose
$z(n)\in
\mathbb{Z}$
 such that
$\sum_{n=0}^k z(n) a(n)>k$.
Then the series
$\sum_{n=0}^\infty z(n) a(n)$ diverges, thereby proving that
$A$ is not \ap\ in $\mathbb{R}$.

To prove the ``only if'' part, assume that $\Cp{X}{\mathbb{R}}$ is
not \IANS\, and choose an infinite \ap\ set
$A\subseteq\Cp{X}{\mathbb{R}}$. Let $\{x_i:i\in\mathbb{N}\}\subseteq
X$ and $\{f_i:i\in\mathbb{N}\}\subseteq A$ be as in the conclusion
of Lemma \ref{selecting:points} (with $0$ instead of $e$ due to the
additive notation). By induction on $n\in\mathbb{N}$ select
$z(n)\in\mathbb{
Z}$ so that
\begin{equation}
\label{eq:6.13:19} \sum_{i=0}^n z(i)f_i(x_n)>n.
\end{equation}
Since $A$ is \ap, there exists $f\in \Cp{X}{\mathbb{R}}$ such that
$f=\lim_{k\to\infty} \sum_{i=0}^k z(i)f_i$. From (\ref{eq2}) and
\eqref{eq:6.13:19} we get $f(x_n)=\sum_{i=0}^n
 z(i)f_i(x_n)>n$ for every $n\in\mathbb{N}$. Thus, the function $f$
is unbounded on $X$, and so $X$ is not pseudocompact.
\end{proof}

Since $\mathbb{R}$-equivalence coincides with $l$-equivalence, from
Theorem \ref{Rpseudo}
we obtain the
following
well-known result of Arhangel'ski\u{\i}.
\begin{corollary}
\label{R:equivalence:preserves:pseudocompactness}
{\rm \cite{Arkhangel}}
$l$-equivalence preserves pseudocompactness.
\end{corollary}

\section{Pseudocompactness of $X$ and \IANS\ property of $C_p(X,G)$}
\label{pseu}

The main goal of this section is to generalize Theorem \ref{Rpseudo}
by replacing the real line $\mathbb{R}$ in it with an arbitrary
NSS group $G$ (provided that $X$ is $G$-regular), see Theorem
\ref{main}.

\begin{lemma}\label{cc}
If $X$ is a countably compact space and  $G$ is an NSS group,
then
$\Cp{X}{G}$ has
the
\IANS\ property.
\end{lemma}
\begin{proof}
Theorem \ref{NSSjeTAP} yields that $G$ is \IANS.
Suppose
that  $A$ is an infinite \ap\  set in
$\Cp{X}{G}$.
Let $\{x_i:i\in\mathbb{N}\}\subseteq X$
and
$\{f_i:i\in\mathbb{N}\}\subseteq A$ be as in the conclusion of
Lemma \ref{selecting:points}.
Since $X$ is countably compact, the set $\{x_i:i\in\mathbb{N}\}$
has a cluster point $x\in X$.
By Lemma
\ref{konecne} the set $J=\{j\in\mathbb{N}: f_j(x)\not=e\}$ is finite.
Let $j=\max J$. After deleting the first $(j+1)$-many $f_i$'s and renumbering, we can assume, without loss of the generality, that
\begin{equation}
\label{eq:f_i}
f_i(x)=e \mbox{ for every }i\in\mathbb{N}.
\end{equation}

Since $G$ is NSS, there exists an open neighborhood $U$ of $e$ in
$G$
as in
item (ii) of
Lemma
 \ref{charNSS}. By recursion on $n\in\mathbb{N}$ we will choose
$z_n\in\mathbb{Z}$ such that
\begin{equation}
\label{**:n}
\tag{**$_n$}
\prod_{i=0}^{n}f_i(x_n)^{z_i} \not\in U.
\end{equation}
Indeed, applying item (ii) of
Lemma
 \ref{charNSS} to $g=f_0(x_0)\not=e$ and $f=h=e$, we can select
$z_0\in \mathbb{Z}$ satisfying (**$_0$). Let
$n\in\mathbb{N}\setminus\{0\}$,
 and suppose
that $z_i\in\mathbb{Z}$ satisfying  (**$_i$)
have already been selected for each $i\in\mathbb{N}$ with $i<n$.
Applying item (ii) of
Lemma
 \ref{charNSS} to $f=e$, $g=f_n(x_n)$ and
$h=\prod_{i=0}^{n-1}f_i(x_n)^{z_i}$, we can find $z_n\in\mathbb{Z}$
satisfying (**$_n$).

Since
$A$ is an \ap\ subset of $\Cp{X}{G}$,
there exists $f\in \Cp{X}{G}$ such that
\begin{equation}
\label{eq:f}
f=\prod_{i=0}^\infty f_i^{z_i}.
\end{equation}
From this and
(\ref{eq2})
we conclude that
$$
f(x_n)=\lim_{k\to\infty} \prod_{i=0}^{k}f_i(x_n)^{z_i}=
\prod_{i=0}^nf_i(x_n)^{z_i}
\mbox{ for every }
n\in\mathbb{N}.
$$
 Combining this with
(**$_n$), we get $f(x_n)\not\in U$ for every $n\in\mathbb{N}$.
Since $x$ is a cluster point of the set
$\{x_n:n\in\mathbb{N}\}$ and $f$ is continuous,
$f(x)$ must be a cluster point of the set
$\{f(x_n):n\in\mathbb{N}\}$, which yields
$f(x)\not\in U$.
On the other hand,
from \eqref{eq:f_i} and \eqref{eq:f} we should have
$$
f(x)=\lim_{k\to\infty} \prod_{i=0}^{k}f_i(x)^{z_i}=e\in U,
$$
a contradiction. This proves that all \ap\ subsets of $\Cp{X}{G}$ are finite.
\end{proof}

\begin{lemma}\label{charofpsc}
If $X$ is a pseudocompact
space and $G$ is a metrizable NSS group,
then $\Cp{X}{G}$ is \IANS.
\end{lemma}
\begin{proof}
Assume that $\Cp{X}{G}$ is not \IANS,
and let $\mathcal{F}$ be an infinite \ap\ subset of $\Cp{X}{G}$.
Lemma \ref{subset:of:productive:is:prodictive}(i) allows us to assume,
without loss of generality, that $\mathcal{F}$ is countable, so
we can fix a faithful enumeration
$\mathcal{F}=\{f_n:n\in\mathbb{N}\}$ of $\mathcal{F}$.
Let $h=\bigtriangleup_{n\in\mathbb{N}}f_n:X\rightarrow G^{\mathbb{N}}$
be the diagonal product.
Since each $f_n$ is
continuous, so is $h$. Since $X$ is pseudocompact, $Y=h(X)$ is
pseudocompact as well. Being a subspace of the metrizable space
$G^{\mathbb{N}}$, $Y$ is metrizable. It follows that $Y$ is compact.

For $n\in \mathbb{N}$ let $p_n: G^{\mathbb{N}}\to G$ be the
projection on $n$th coordinate (defined by $p_n(\phi)=\phi(n)$ for
$\phi\in G^{\mathbb{N}}$). Since $p_n$ is continuous,
$g_n=p_n\upharpoonright_Y\in\Cp{Y}{G}$. Clearly, $f_n=g_n \circ
h$.
If $m,n\in\mathbb{N}$ and $m\not=n$, then
$f_m\not=f_n$,
 and so
$g_m(h(x))=g_m \circ
h(x)=f_m(x)\not=f_n(x)=g_n \circ h(x)=g_n(h(x))$ for some $x\in X$,
which yields $g_m\not=g_n$. Therefore, the family
$\mathcal{G}=\{g_n:n\in\mathbb{N}\}\subseteq \Cp{Y}{G}$
is
faithfully indexed (in particular, infinite).

Let us show
that $\mathcal{G}$
 is
\ap\ in $\Cp{Y}{G}$, in contradiction with Lemma \ref{cc}. Let
$s:\mathbb{N}\to\mathbb{N}$ be an injection and
$z:\mathbb{N}\to\mathbb{Z}$ a map. Since $\mathcal{F}$ is \ap,
there exists $f\in\Cp{X}{G}$ such that $f=\prod_{n=0}^\infty
f_{s(n)}^{z(n)}$.

Assume that $x,x'\in X$ and $h(x)=h(x')$. For each $n\in \mathbb{N}$
we have $f_n=g_n \circ h$, which yields $f_n(x)=f_n(x')$. Therefore,
$$
f(x)=\lim_{k\to \infty} \prod_{n=0}^k f_{s(n)}^{z(n)}(x)= \lim_{k\to
\infty} \prod_{n=0}^k f_{s(n)}^{z(n)}(x')=f(x').
$$
It follows that there exists a unique function $g:Y\to G$ such that
$f=g\circ h$.
Since $X$ is pseudocompact, $G$ is metrizable and $f\in\Cp{X}{G}$,
from \cite[Theorem 7]{R-quotient} we conclude that $g\in \Cp{Y}{G}$.

Let $y\in Y$ be arbitrary. Choose $x\in X$ such that $y=h(x)$. Then
\begin{align*}
g(y)=g(h(x))= f(x)&= \lim_{k\to \infty} \prod_{n=0}^k f_{s(n)}^{z(n)}(x)=
\lim_{k\to \infty} \prod_{n=0}^k (g_{s(n)} \circ h)^{z(n)}(x)\\
&= \lim_{k\to
\infty} \prod_{n=0}^k \left(g_{s(n)} (h(x))\right)^{z(n)}= \lim_{k\to \infty}
\prod_{n=0}^k g_{s(n)} ^{z(n)}(y).
\end{align*}
Therefore,
$g=\prod_{n=0}^\infty g_{s(n)}^{z(n)}$.
Thus, $\mathcal{G}$ is \ap\ in $\Cp{Y}{G}$.
\end{proof}

\begin{theorem}\label{mainlemma}
If $X$ is a pseudocompact space and $G$ is an
NSS group,
then
$\Cp{X}{G}$ is \IANS.
\end{theorem}
\begin{proof}
Assume that $\Cp{X}{G}$ is not \IANS, and let $A$ be an infinite \ap\ subset of $\Cp{X}{G}$.
Lemma \ref{subset:of:productive:is:prodictive}(i) allows us to assume,
without loss of generality, that $A$ is countable, so
we can fix a faithful enumeration
$A=\{f_n:n\in\mathbb{N}\}$.
For each $n\in\N$, the subset $f_n(X)$ of $G$ is
pseudocompact, being a continuous image of the pseudocompact space
$X$.
The smallest closed subgroup $K$
of $G$ containing $\bigcup \{f_n(X):n\in\N\}$ must be $\omega$-bounded.
(Recall that, according to Guran \cite{Guran},
 a topological group $G$ is called {\em $\omega$-bounded\/}
if, for any open set $U\subseteq G$, there exists a countable set
$S\subseteq G$
such that $SU=\{su:s\in S, u\in U\}=G$.)
Note that
$A\subseteq \Cp{X}{K}$ and
$\Cp{X}{K}$ is a closed subgroup of $\Cp{X}{G}$, so
$A$ is \ap\ in $\Cp{X}{K}$ by Lemma
\ref{basic:properties:of:ap}(ii).
Being a subgroup of an NSS group $G$, $K$ itself is an
NSS group. Therefore, without loss of generality, we may (and will)
assume that $K=G$, i.e., $G$ itself is $\omega$-bounded. Thus, there
exists a family $\{G_\beta:\beta\in B\}$ consisting of separable
metric groups $G_\beta$ such that $G$ is a subgroup of
$\prod_{\beta\in B}G_\beta$; see \cite{Guran}.

Let $n,m\in \N$ and $n< m$. Since $f_n\neq f_m$, there is
$x_{n,m}\in X$ such that $g_{n,m}=f_n(x_{n,m})\neq
f_m(x_{n,m})=g_{m,n}$. So we can pick some $\beta_{n,m}\in B$ such
that $g_{n,m}(\beta_{n,m})\not=g_{m,n}(\beta_{n,m})$.
Since
$G$ is an NSS group, we can find $k\in\mathbb{N}$,
$\beta_1,\dots,\beta_k\in B$
and an open neighborhood $U_i$ of
the identity in $G_{\beta_i}$ for $i\le k$, such that the open
neighborhood
$$
U=G\cap\left\{g\in\prod_{\beta\in B}G_\beta: g(\beta_i)\in U_i\mbox{ for
}i=1,\dots,k\right\}
$$
of the identity of $G$ contains no nontrivial subgroup of $G$.
Define
$$
C=\{\beta_{n,m}:n,m\in \N,n< m\}\cup\{\beta_1,\dots,\beta_k\}
$$
and consider the
projection $q:\prod_{\beta\in B}G_\beta\rightarrow\prod_{\beta\in
C}G_\beta$. Let $H=q(G)$ and $\phi=q\upharpoonright_G$. As a
subspace of a countable product of metrizable spaces, $H$ is
metrizable. Moreover,
$\phi(U)$ is an
open neighborhood of the identity of $H$ which contains no
nontrivial subgroup of $H$. Hence, $H$ is NSS.

Let $\Phi:\Cp{X}{G}\rightarrow\Cp{X}{H}$ be the continuous
homomorphism
defined
by $\Phi(f)=\phi\circ f$
for $f\in\Cp{X}{G}$.
Then $\Phi(A)$
is an \ap\ subset of $\Cp{X}{H}$ by Lemma \ref{a.p.ispreserved}(ii).

Let $m,n\in\mathbb{N}$ and $n<m$.
Since
$\beta_{n,m}\in C$ and
$g_{n,m}(\beta_{n,m})\not=g_{m,n}(\beta_{n,m})$, we have
$q(g_{n,m}) \not= q(g_{m,n})$, which yields
$$
\Phi(f_n)(x_{n,m})=\phi( f_n(x_{n,m}))=\phi(g_{n,m})=q(g_{n,m})
\not=$$$$ q(g_{m,n})=\phi(g_{m,n})=\phi(
f_m(x_{n,m}))=\Phi(f_m)(x_{n,m}).
$$
Hence $\Phi(f_n)\not=\Phi(f_m)$.
This shows that $\Phi(A)=\{\Phi(f_n):n\in\mathbb{N}\}$
is a faithfully indexed (and thus infinite) set.
It follows that $\Cp{X}{H}$ is not \IANS, in contradiction with Lemma \ref{charofpsc}.
\end{proof}
\begin{lemma}\label{pseudo}
Let $G$ be a
topological group and $X$ a $G$-regular
space which is not pseudocompact. Then $C_p(X,G)$ contains a
subgroup topologically isomorphic to the product
$H=\prod_{i\in\mathbb{N}} H_i$ of nontrivial topological groups
$H_i$.
\end{lemma}
\begin{proof}
Since $X$ is not pseudocompact, there exists a discrete family
$\mathscr{U}=\{U_i:i\in\mathbb{N}\}$ consisting of non-empty open
subsets of $X$. For each $i\in\mathbb{N}$ choose $x_i\in U_i$ and
use $G$-regularity of $X$ to fix $f_i\in C_p(X,G)$ such that
$f_i(X\setminus U_i)\subseteq\{e\}$ and $f_i(x_i)\neq e$. Let $H_i$ be
the cyclic subgroup of $C_p(X,G)$ generated by $f_i$ (equipped with
the subspace topology inherited from $C_p(X,G)$). Clearly, $H_i$ is
nontrivial. Let $H=\prod_{i\in\mathbb{N}} H_i$. Since the family
$\mathscr{U}$ is discrete, for each
$h\in H$
the infinite product
$\theta(h)=\prod_{i=0}^\infty h(i)=\lim_{k\to\infty} \prod_{i=0}^k
h(i)$
 is well-defined and $\theta(h)\in C_p(X,G)$. A straightforward
verification of the fact that $\theta: H\to C_p(X,G)$ is a
topological isomorphism between $H$ and $\theta(H)$ is left to the
reader.
\end{proof}

\begin{theorem}\label{main}
For an
NSS group $G$ and a $G$-regular space $X$, the
following conditions are equivalent:
\begin{itemize}
\item[(i)] $X$ is pseudocompact;
\item[(ii)]$\Cp{X}{G}$ is \IANS;
\item[(iii)]$\Cp{X}{G}$ does not contain a subgroup which is topologically isomorphic to an
infinite product of nontrivial
 topological
 groups.
\end{itemize}
\end{theorem}
\begin{proof}
(i)$\Rightarrow$(ii) follows from Theorem
\ref{mainlemma}.

(ii)$\Rightarrow$(iii)
is Proposition \ref{h}(iii).

(iii)$\Rightarrow$(i) follows from
Lemma \ref{pseudo}.
\end{proof}

\begin{corollary}\label{opseudocomp}
Let $G$ be an NSS  group. Then $G$-equivalence preserves
pseudocompactness within the class of $G$-regular spaces.
\end{corollary}

The proof of the following lemma is straightforward.
\begin{lemma}\label{Dima}
Let $X$ and $Y$ be spaces and $H$ a topological group. For
$f\in C_p(X\times Y, H)$
and $x\in X$ define $f_x\in \Cp{Y}{H}$ by
$f_x(y)=f(x,y)$ for every $y\in Y$.
Consider the map $\theta: C_p(X\times Y, H)\to C_p(X,C_p(Y,H))$
which assigns to every $f\in C_p(X\times Y, H)$ the function
$\theta(f)\in C_p(X,C_p(Y,H))$ defined by $\theta(f)(x)=f_x$ for each $x\in X$.
Then $\theta$
is a topological isomorphism
between $C_p(X\times Y, H)$ and  $\theta(C_p(X\times Y, H))$.
\end{lemma}

It follows from Propositions \ref{closed:copy:of:G:in:Cp}
and \ref{h}(i)
that
$G$ must be \IANS\ whenever $\Cp{X}{G}$ is \IANS.
Our next theorem shows that the \IANS\ property of $G$ is not sufficient to ensure that
$\Cp{X}{G}$ is \IANS.

\begin{theorem}
\label{separately:continuous:example}
There exist a
precompact \IANS\ group $G$
and
 a countably compact $G^\star$-regular space $X$ such that $\Cp{X}{G}$ is not \IANS.
\end{theorem}
\begin{proof}
Let $X$ be a countably compact Tychonoff space
and $Y$ a pseudocompact Tychonoff space such that $X\times Y$ is not
pseudocompact (see, for example, \cite[Example 3.10.19]{Engelking}).

By Proposition \ref{Gregualrity}, $Y$ is  $\mathbb{T}$-regular (in
fact, even $\mathbb{T}^{\star\star}$-regular). Since $Y$ is
pseudocompact and $\mathbb{T}$ is NSS, Theorem \ref{main} yields
that $G=C_p(Y,\mathbb{T})$ is a \IANS\ group. Since $G$ is a
subgroup of the compact group $\mathbb{T}^Y$, $G$ is precompact. By
Proposition \ref{closed:copy:of:G:in:Cp}, $G$ contains a subgroup
topologically isomorphic to $\mathbb{T}$, so $X$ is
$G^\star$-regular by Proposition \ref{Gregualrity}(ii).

By Proposition \ref{Gregualrity}, the space $X\times Y$ is
$\mathbb{T}$-regular (in fact, even $\mathbb{T}^{\star\star}$-regular).
Since $X\times Y$ is not pseudocompact, it follows from Lemmas \ref{pseudo}
and \ref{algebraically:productive:subsets:in:products}(i)
that
$\Cp{X\times Y}{\mathbb{T}}$ contains an infinite set $A$ that is \ap\  in $\Cp{X\times Y}{\mathbb{T}}$.
According to Lemma \ref{Dima}, $\Cp{X\times Y}{\mathbb{T}}$ is
topologically isomorphic to a subgroup of
$\Cp{X}{\Cp{Y}{\mathbb{T}}}=\Cp{X}{G}$.
Applying Lemma \ref{basic:properties:of:ap}(i), we conclude that $A$ is \ap\
in $\Cp{X}{G}$. Since $A$ is infinite, it follows that $\Cp{X}{G}$ is not
\IANS.
\end{proof}

Theorem
\ref{separately:continuous:example}
demonstrates that
the conclusion of
Theorem \ref{main} is no longer valid if we
replace the NSS property of $G$ in its assumption by the weaker \IANS\ property.

\begin{proposition}
\label{NoNSS}
\label{Cp:TAP:not:NSS}
Let $G$ be a
topological group and $X$ an infinite
$G$-regular space.
\begin{itemize}
\item[(i)]
$\Cp{X}{G}$ is not NSS.
\item[(ii)]
If, in addition, $X$ is pseudocompact and $G$ is NSS,
 then
$\Cp{X}{G}$ is \IANS\ but not NSS.
\end{itemize}
\end{proposition}
\begin{proof}
(i)
Let $U$ be any neighborhood of the identity in $\Cp{X}{G}$. Then
there exist an integer $n\in\mathbb{N}\setminus\{0\}$, points $x_1,\ldots,x_n\in X$ and an open set $V\subseteq G$ with $e\in
V$, such that
$H=\{f\in\Cp{X}{G}: f(x_i)
=e$ for  $i=1,\ldots,n\}
\subseteq U$.
Since $X$ is infinite and $G$-regular, $H$ is a nontrivial subgroup
of $\Cp{X}{G}$. It follows that $\Cp{X}{G}$ is not NSS.

(ii) follows from (i) and Theorem \ref{mainlemma}.
\end{proof}

\section{Compactness-like properties and
$G$-equivalence}\label{seccomp}

For a space $X$, $nw(X)$ denotes the network weight of $X$, $t(X)$
stays for the tightness of $X$, $l(X)$ denotes the Lindel\"{o}f
number of $X$, and $l^*(X)=sup\{l(X^n):n\in\mathbb{N}\}$.

We start with a ``$G$-analogue'' of the well-known cardinal equality from the $C_p$-theory.
\begin{lemma}
\label{network:lemma}
If $G$ is a separable metric group, then $nw(X)=nw(\Cp{X}{G})$
for every $G$-regular space $X$.
\end{lemma}
\begin{proof}
Let $\mathcal{N}$ be a network of $X$ such that
$|\mathcal{N}|\le nw(X)$, and let $\mathcal{B}$ be a countable base of
$G$.  For $N\in\mathcal{N}$ and $U\in\mathcal{B}$
define $W(N,U)=\{f\in\Cp{X}{Y}:f(N)\subseteq U\}$. Then the family
consisting of finite intersections of the members of the family
$\mathcal{W}=\{W(N,U):N\in\mathcal{N}, U\in\mathcal{B}\}$ is a
network of $\Cp{X}{G}$ satisfying
$|\mathcal{W}|\le|\mathcal{N}|\cdot\omega\le nw(X)$. This proves the
inequality $nw(\Cp{X}{G})\leq nw(X)$ for an {\em arbitrary\/} (not
necessarily $G$-regular) space $X$. In particular,
$nw(\Cp{\Cp{X}{G}}{G})\leq nw(\Cp{X}{G})$.

For
$x\in X$, let $\pi_x:\Cp{X}{G}\to G$ be the projection defined by
$\pi_x(f)=f(x)$ for $f\in\Cp{X}{G}$. A straightforward check, using $G$-regularity of $X$, shows that the map
$\phi:X\to\Cp{\Cp{X}{G}}{G}$ given by $\phi(x)=\pi_x$ for $x\in X$, is a
homeomorphism. Therefore, $nw(X)=nw(\phi(X))\le nw(\Cp{\Cp{X}{G}}{G})\leq nw(\Cp{X}{G})\le nw(X)$.
\end{proof}

\begin{corollary}\label{networkweight:cor}
Let $G$ be a separable metric group. Then $G$-equivalence preserves
the
network weight within the class of $G$-regular spaces.
\end{corollary}

Our next lemma is a ``$G$-analogue'' of the well-known theorem of
Pytkeev from the $C_p$-theory; see \cite[Theorem 1]{Pytkeev}. Its
proof essentially follows the original proof, with necessary
adaptations to take into account the $G^\star$-regularity condition.
\begin{lemma}
\label{lemma:8.2}
If
$G$ is a topological group and $X$ is a $G^\star$-regular space,
then $l^*(X)\le t(\Cp{X}{G})$.
\end{lemma}
\begin{proof}
To start with, we claim that one can find  $g\in G\setminus\{e\}$
such that, for every open subset $U$ of $X$ and each non-empty
finite set $K\subseteq U$, there exists $f_{K,U}\in\Cp{X}{G}$
satisfying $f_{K,U}(K)\subseteq \{e\}$ and $f_{K,U}(X\setminus U)\subseteq
\{g^{-1}\}$. Indeed, let $g\in G\setminus\{e\}$ be the element
witnessing $G^\star$-regularity of $X$. Let $U$ be an open subset of
$X$ and $K\not=\emptyset$ a finite subset of $U$. Let
$K=\{x_0,\dots,x_k\}$ be a faithful enumeration of $K$, and let
$U_0,\dots,U_k$ be pairwise disjoint open subsets of $U$ with
$x_i\in U_i$ for $i\le k$. For every $i\le k$ we can choose
$\varphi_i\in \Cp{X}{G}$ such that $\varphi_i(x_i)=g$ and
$\varphi_i(X\setminus U_i)\subseteq\{e\}$. Then the function
$f_{K,U}\in\Cp{X}{G}$ defined by
$f_{K,U}(x)=g^{-1}\cdot\prod_{i=0}^k \varphi_i(x)$ for $x\in X$, is
as required.

Given two families $\mathscr{A}$ and $\mathscr{B}$, we write
$\mathscr{A}\prec \mathscr{B}$ provided that, for every $A\in \mathscr{A}$,
there exists $B\in\mathscr{B}$ such that $A\subseteq B$.

Fix $n\in\mathbb{N}\setminus\{0\}$ and an open cover $\mathscr{V}$ of
$X^n$.
Let $\mathbf{U}$ denote the set of all finite families $\mathscr{U}$ of open subsets of $X$
satisfying $\Pi(\mathscr{U})\prec\mathscr{V}$, where
$\Pi(\mathscr{U})=\{U_1\times\dots\times U_n:(U_1,\dots,U_n)\in\mathscr{U}^n\}$.
For every $\mathscr{U}\in\mathbf{U}$ choose a finite subfamily
$\mathscr{V}_{\mathscr{U}}$ of $\mathscr{V}$ such that
$\Pi(\mathscr{U})\prec\mathscr{V}_{\mathscr{U}}$.
Let
\begin{equation}
\label{eq:19}
F=\left\{f\in \Cp{X}{G}:
f\left(X\setminus\bigcup\mathscr{U}\right)\subseteq \{g^{-1}\}
\mbox{ for some }
\mathscr{U}\in \mathbf{U}
\right\}.
\end{equation}

We claim that $\mathbf{1}\in \overline{F}$, where
$\mathbf{1}\in\Cp{X}{G}$ is defined by $\mathbf{1}(x)=e$ for all
$x\in X$. Indeed, let $K$ be a non-empty finite subset of $X$ and
$O$ an arbitrary open neighborhood of $e$ in $G$. Since
$\mathscr{V}$ is an open cover of $X^n$, there exists some
$\mathscr{U}\in\mathbf{U}$ with $K\subseteq U$, where
$U=\bigcup\mathscr{U}$. Now $f_{K,U}\in F\cap\bigcap_{x\in K}
W(x,O)\not=\emptyset$.

Since $\mathbf{1}\in \overline{F}$,
we can choose $F^*\subseteq F$ such that
$|F^*|\le t(\Cp{X}{G})$
and $\mathbf{1}\in \overline{F^*}$.
For every $f\in F^*$ use \eqref{eq:19} to select $\mathscr{U}_f\in\mathbf{U}$ such that
\begin{equation}
\label{eq:20}
f\left(X\setminus\bigcup\mathscr{U}_f\right)\subseteq \{g^{-1}\}.
\end{equation}
Define $\mathscr{V}^*=\bigcup_{f\in
F^*}\mathscr{V}_{\mathscr{U}_f}\subseteq \mathscr{V}$. Clearly,
$|\mathscr{V}^*|\le |F^*|\le t(\Cp{X}{G})$. It remains only to prove
that $\mathscr{V}^*$ covers $X^n$. Indeed, let $(x_1,\dots,x_n)\in
X^n$ be arbitrary. Define $K=\{x_1,\dots,x_n\}$ and
$O^*=G\setminus\{g^{-1}\}$. Then $\bigcap_{x\in K} W(x,O^*)$ is an
open neighborhood of $\mathbf{1}$, so we can pick
\begin{equation}
\label{eq:21}
f\in F^*\cap \bigcap_{x\in K} W(x,O^*).
\end{equation}
From \eqref{eq:21} and \eqref{eq:20} we conclude that
$K\subseteq \bigcup\mathscr{U}_f$.
For $i=1,\dots,n$ choose $U_i\in \mathscr{U}_f$ with $x_i\in U_i$.
Then $(x_1,\dots,x_n)\in U_1\times\dots\times U_n\in \Pi(\mathscr{U}_f)$.
Since $\Pi(\mathscr{U}_f)\prec \mathscr{V}_{\mathscr{U}_f}$,
we have
$U_1\times\dots\times U_n\subseteq V$ for some $V\in \mathscr{V}_{\mathscr{U}_f}\subseteq \mathscr{V}^*$.
Hence $(x_1,\dots,x_n)\in\bigcup\mathscr{V}^*$.
\end{proof}

Our next proposition is a ``$G$-analogue'' of the classical theorem of
Arhangel'ski\u{\i}-Pytkeev from the $C_p$-theory.
\begin{proposition}\label{tightness}
If $G$ is a metric group and $X$ is a $G^\star$-regular space, then
$l^*(X)=t(\Cp{X}{G})$.
\end{proposition}
\begin{proof}
The inequality $l^*(X)\le t(\Cp{X}{G})$
was proved in Lemma \ref{lemma:8.2}.
The converse inequality can be proved by a straightforward modification of
the proof of \cite[Theorem II.1.1]{A}.
\end{proof}

\begin{corollary}
\label{corollary:8.4}
Let $G$ be a metric group. If  $G^\star$-regular spaces $X$ and $Y$ are $G$-equivalent, then $l^*(X)= l^*(Y)$.
\end{corollary}

Theorem \ref{main} and Proposition \ref{tightness} give the
following
\begin{theorem}\label{char:of:compactness}
Let $G$ be an NSS  metric group. Then
a $G^\star$-regular space
$X$ is compact if and only if $\Cp{X}{G}$ is
a \IANS\ group of countable tightness.
\end{theorem}

\begin{corollary}\label{nademnou}
Let $G$ be an NSS  metric group. Then $G$-equivalence preserves
compactness within the class of $G^\star$-regular spaces.
\end{corollary}

Since a space is compact and metrizable if and
only if it is pseudocompact and has a countable network, combining
Theorem \ref{main} and Lemma \ref{network:lemma}, we get the following
\begin{theorem}\label{char:of:compact:metrizable}
Let $G$ be an NSS separable metric group.
For a $G$-regular space
$X$, the following conditions are equivalent:
\begin{itemize}
 \item[(i)]
  $X$ is compact and metrizable;
 \item[(ii)]
$\Cp{X}{G}$ is a \IANS\ group with a countable network.
\end{itemize}
\end{theorem}
\begin{corollary}
Let $G$ be an NSS separable metric group. Then $G$-equivalence
preserves the property ``to be compact metrizable'' within the class
of $G$-regular spaces.
\end{corollary}

Our next example shows that $(G\times H)$-equivalence need not imply either
$G$-equivalence or $H$-equivalence.

\begin{example}\label{E3}\label{GxH:does:not:imply:G}
\rm{ Let $G=\mathbb{T}^\omega\times\mathbb{R}$ and
$H=\mathbb{T}\times\mathbb{R}^\omega$. Then {\em both $G$ and
$H$-equivalence preserve pseudocompactness and compactness, but
$G\times H$-equivalence does not preserve either of these
properties\/}. By Proposition \ref{Gregualrity}(i), every
space is both $G^{\star\star}$-regular ad
$H^{\star\star}$-regular. Applying either Theorem \ref{thm:4.4}(ii)
or Theorem \ref{thm:4.5}(i), we conclude that $G$-equivalence
implies $\mathbb{R}$-equivalence, and $H$-equivalence implies
$\mathbb{T}$-equivalence. Both $\mathbb{R}$-equivalence and
$\mathbb{T}$-equivalence preserve pseudocompactness (Corollary
\ref{opseudocomp}) and compactness (Corollary \ref{nademnou}). Thus,
both $G$-equivalence and $H$-equivalence preserve pseudocompactness
and compactness as well. On the other hand,  since $G\times
H\cong(\mathbb{R}\times\mathbb{T})^\omega$, Corollary
\ref{Gomegadoesntpreserve} yields that $(G\times H)$-equivalence
preserves neither pseudocompactness, nor compactness. }
\end{example}

In connection with
Corollary \ref{opseudocomp}, our next
example
shows that
neither NSS nor \IANS\ property of $G$ is
necessary for $G$-equivalence to preserve pseudocompactness
(and compactness as well).

\begin{example}
\label{E1}
\begin{itemize}
\item[(i)]
\rm{
{\em For every infinite zero-dimensional
pseudocompact
 space $X$, the group $H=\mathbb{R}\times
\Cp{X}{\mathbb{Z}(2)}$ is \IANS\ but not NSS, and $H$-equivalence
preserves both pseudocompactness and compactness\/}. Indeed,
according to
Theorem \ref{mainlemma},
 $\Cp{X}{\mathbb{Z}(2)}$ is \IANS. Since
$\mathbb{R}$ is \IANS\ too, it follows from Proposition \ref{h}(i)
that $H$ is \IANS\ as well. From Proposition \ref{NoNSS}(i) we get
that $\Cp{X}{\mathbb{Z}(2)}$, and consequently $H$, is not NSS.
Furthermore, applying either Theorem \ref{thm:4.4}(i) or Theorem
\ref{thm:4.5}(ii), we obtain that $H$-equivalence implies
$\mathbb{R}$-equivalence. Since $\mathbb{R}$-equivalence (that is
$l$-equivalence) preserves pseudocompactness (Corollary
\ref{R:equivalence:preserves:pseudocompactness}) and compactness
(Corollary \ref{nademnou}), so does $H$-equivalence. } \label{E2}
\item[(ii)]
\rm{
{\em
The group $G=\mathbb{R}\times\mathbb{Z}(2)^\omega$ is not \IANS\
(and consequently not NSS by Theorem \ref{NSSjeTAP}), but
$G$-equivalence preserves both pseudocompactness and compactness.}
Proposition \ref{h}(iii) guarantees that $G$ is not \IANS.
Applying either Theorem \ref{thm:4.4}(i) or  Theorem \ref{thm:4.5}(ii),
we conclude that $G$-equivalence
implies $\mathbb{R}$-equivalence.
Now we finish the argument as in item (i).
}
\end{itemize}
\end{example}

\begin{remark}
{\rm Examples \ref{E3} and \ref{E1} show that the class of all
(Abelian) topological groups $G$ for which $G$-equivalence preserves
compactness is not finitely productive and is larger then that of
NSS, metrizable groups. }
\end{remark}

Now we turn to a particular version of Problem \ref{problem:3}
by considering the class $\mathbf{PSC}$ of all topological groups $G$ for which
$G$-equivalence preserves pseudocompactness.

By Corollary \ref{opseudocomp},
$\mathbb{Z}(2)$-equivalence preserves
 pseudocompactness within the class
of  $\mathbb{Z}(2)$-regular spaces, and yet $\mathbb{Z}(2)\not\in
\mathbf{PSC}$, by Example \ref{negative:example}. Therefore, it is
reasonable to investigate whether $G$-equivalence preserves
pseudocompactness only within the class of $G$-regular spaces. On
the other hand, if a space $X$ is $G$-regular for every topological
group $G$, then $X$ must be zero-dimensional. One possible way to
avoid such a restriction on $X$ is to require our groups $G$ to be
the elements of the class $\mathbf{I}$ of all topological groups
that contain a homeomorphic copy of the closed unit interval $[0,1]$
as a subspace. Indeed, by Proposition \ref{Gregualrity}, this would
make the condition of $G$-regularity automatically satisfied for
every space. Therefore, one may expect that the subclass
$\mathbf{PSC}\cap\mathbf{I}$ of the class $\mathbf{PSC}$ should have
especially nice properties. Let us summarize what we know about the
properties of this class.

\begin{proposition}
Denote by $\mathbf{NSS}$ the class of all
NSS groups and by $\mathbf{TAP}$ the class of all \IANS\ groups.
Then:
\begin{itemize}
\item[(i)]
$\mathbf{NSS}\cap\mathbf{I}\subseteq\mathbf{PSC}\cap\mathbf{I}$
(Corollary \ref{opseudocomp});
\item[(ii)]
$\mathbf{PSC}\cap\mathbf{I}\neq\mathbf{I}$ (Corollary
\ref{Gomegadoesntpreserve});
\item[(iii)]
$(\mathbf{PSC}\cap\mathbf{I})\setminus\mathbf{NSS}\not=\emptyset$
(Example \ref{E1}(i));
\item[(iv)]
$(\mathbf{PSC}\cap\mathbf{I})\setminus\mathbf{TAP}\not=\emptyset$
(Example \ref{E2}(ii));
\item[(v)]
both $\mathbf{PSC}$ and $\mathbf{PSC}\cap\mathbf{I}$ are closed under
taking finite powers (Corollary \ref{productOK});
\item[(vi)]
both $\mathbf{PSC}$ and $\mathbf{PSC}\cap\mathbf{I}$ are not closed under taking finite
products (Example \ref{E3}).
\end{itemize}
\end{proposition}

\section{(Dis)connectedness and $G$-equivalence}\label{secdis}

Recall that a topological space $X$ is {\em totally disconnected\/}
if every quasi-component of $X$ is a singleton, or equivalently,
if for
every
pair $x,y$ of distinct
points
of $X$
there exists a clopen
set $F\subseteq X$ such that $x\in F\not\ni y$.

\begin{theorem}\label{totdiscon}
Let $G$ be a topological group with the dense and totally disconnected
torsion part $tor(G)$. Then a $G$-regular space $X$ is totally
disconnected if and only if $tor(\Cp{X}{G})$ is dense in
$\Cp{X}{G}$.
\end{theorem}
\begin{proof}
Suppose that $X$ is totally disconnected.
Let $O$ be a non-empty open subset of $\Cp{X}{G}$.
Choose
$f\in O$.
 Then there exist $n\in\mathbb{N}\setminus\{0\}$, pairwise
distinct elements $x_1,\dots,x_n$ of $X$ and non-empty open subsets
$U_1,\dots,U_n$ of $G$ such that $f\in\bigcap_{i=1}^n
W(x_i,U_i)\subseteq O$. As $tor(G)$ is dense in $G$, for every
$i=1,\ldots , n$ we can choose $t_{i}\in tor(G)\cap U_i$. Since $X$
is totally disconnected, there exists a disjoint partition
$X=\bigcup_{i=1}^n F_i$ of $X$
into clopen subsets $F_i$
 such that $x_{i}\in F_i$ for
$i=1,\ldots , n$. Define $h\in\Cp{X}{G}$ by letting $h(x)=t_i$
whenever $x\in F_i$. Clearly,
 $h\in tor(\Cp{X}{G})$
and $h\in \bigcap_{i=1}^n W(x_i,U_i)\subseteq O$.
Therefore, $O\cap tor(\Cp{X}{G})\not=\emptyset$.
This proves that $tor(\Cp{X}{G})$ is dense in $\Cp{X}{G}$.

To prove the reverse implication, assume that $tor(\Cp{X}{G})$ is dense in
$\Cp{X}{G}$.
Suppose that $x,y\in X$ and $x\not=y$.
Since $X$ is $G$-regular,
 there exists $f\in\Cp{X}{G}$ with
$f(x)\neq f(y)$. Since $tor(\Cp{X}{G})$ is dense in $\Cp{X}{G}$, we
may assume that $f\in tor(\Cp{X}{G})$, and thus $f\in\Cp{X}{tor(G)}$.
Since $tor(G)$ is totally disconnected, there exists a clopen
set $W\subseteq tor(G)$ with $f(x)\in W\not\ni f(y)$. Consequently,
$f^{-1}(W)$ is  a clopen subset of $X$ satisfying $x\in f^{-1}(W)\not\ni y$. This shows that $X$ is totally disconnected.
\end{proof}
\begin{corollary}
Let $G$ be a topological group with the dense and totally disconnected
torsion part $tor(G)$. Then $G$-equivalence preserves total
disconnectedness within the class of $G$-regular spaces.
\end{corollary}

The proof of the following lemma is straightforward.
\begin{lemma}
\label{9.3}
If $X=\bigcup_{i\in I} X_i$
is a decomposition of a space $X$ into a disjoint union of its
non-empty clopen subsets $X_i$, then $\Cp{X}{G}\cong\prod_{i\in I}
\Cp{X_i}{G}$ for every topological group $G$.
\end{lemma}

Recall that the  {\em order\/} of an element $g$ of a group $G$ is
defined to be the smallest $n\in\mathbb{N}\setminus\{0\}$ satisfying
$g^n=e$ (if such $n$ exists).

\begin{proposition}\label{souvislost}
Suppose that
a topological group $G$ has
exactly $m$ elements of order $p$, for a suitable integer
$m\in\mathbb{N}\setminus\{0\}$  and some prime number $p$. Let $k\in
\mathbb{N}\setminus\{0\}$. Then a space $X$ has precisely $k$
connected components if and only if $\Cp{X}{G}$ has exactly
$(m+1)^{k}-1$ elements of order $p$.
\end{proposition}
\begin{proof}
To prove the ``only if'' part,
it suffices to realize that an element $f\in\Cp{X}{G}$ of order
$p$ must
be constant on every connected component of $X$. The rest is a
simple computation.

To prove the ``if'' part, it suffices to show that $X$ has finitely
many connected components, since the (finite) number $m$ of
connected components of $X$ is uniquely determined by $(m+1)^{k}-1$
(this follows from  the ``only if'' part of our proof). Assume that $X$ has
infinitely many pairwise distinct connected components. Then for
every natural number $n\ge 2$ one can find a decomposition
$X=\bigcup_{i=1}^n X_i$ of $X$ into non-empty clopen subsets $X_i$
of $X$, and Lemma \ref{9.3} yields that $\Cp{X}{G}\cong\prod_{i=1}^n
\Cp{X_i}{G}$. Since each $\Cp{X_i}{G}$ has at least one element of
order $p$ by Proposition \ref{closed:copy:of:G:in:Cp},
it follows that $\Cp{X}{G}$ must have at least $n$
elements of order $p$. Since $n$ was chosen arbitrarily, this
contradicts our assumption that $\Cp{X}{G}$ has finitely many
elements of order $p$.
\end{proof}

\begin{corollary}
\label{cor:souvislost} Let $G$ be a topological group that contains
at least one, but finitely many, elements of prime order $p$ (for a
suitable $p$). Then $G$-equivalence preserves the finite number of
connected components.
\end{corollary}

Corollary  \ref{cor:souvislost}
is a
slight
generalization of the result of
Tkachuk,
 who proved in \cite{Tkachuk} that, for every $n\in\mathbb{N}\setminus\{0\}$,
the property
``to consist of $n$-many connected components''
is
preserved by $\mathbb{Z}(2)$-equivalence (and thus, by
$M$-equivalence as well).

For a topological group $G$, we
denote by
$c_0(G)$ the {\em pathwise connected component\/} of $G$ (that is,
 the
union of all pathwise connected subsets of $G$ containing $e_G$),
and we use $c(G)$ for denoting the {\em connected component} of $G$
(that is, the union of all connected subsets of $G$ containing
$e_G$). It is known that $c(G)$ is a closed normal subgroup of $G$.

For certain topological groups $G$, we can characterize the finite
number of connected components of a space $X$ in terms of a purely
topological property of $\Cp{X}{G}$. Recall that, for a natural
number $n\ge 1$,
 a subgroup $H$ of $G$ {\em has index\/ $n$ in $G$\/} provided that
$|G/H|=n$.
\begin{proposition}
\label{souvislost2} Let $m\geq1$ and $n\geq2$ be natural numbers.
Assume that $G$ is a topological group
such that $c_0(G)=c(G)$ has index $n$ in $G$. Then a space $X$ has
precisely $m$ connected components if and only if $c(\Cp{X}{G})$ has
index  $n^m$ in $\Cp{X}{G}$.
\end{proposition}
\begin{proof}
We will start with the ``only if'' part. Let $G=\bigcup_{j=1}^nC_j$,
where $C_1,\ldots,C_n$ are pairwise disjoint translates of $c_0(G)$.
Let $X=\bigcup_{i=1}^mX_i$, where $X_1,\ldots,X_m$ are pairwise
distinct connected components of $X$. Fix $i=1,\dots,m$. Since an
image of a connected space is connected,
$\Cp{X_i}{G}=\bigcup_{j=1}^n\Cp{X_i}{C_j}$. Since each $C_j$ is a
translate of $c_0(G)$, the spaces $\Cp{X_i}{C_j}$ and
$\Cp{X_i}{c_0(G)}$ are homeomorphic. Using the same argumentation as
in the proof of Theorem \ref{thm:4.4} Case 1 we conclude that the
latter group is connected. Hence, each $\Cp{X_i}{C_j}$ is a
connected component of $\Cp{X_i}{G}$. Therefore, the group
$\Cp{X_i}{G}$ has precisely $n$ connected components. Since each of
the finitely many components $X_i$ of $X$ is clopen in $X$, we have
$\Cp{X}{G}\cong\prod_{i=1}^m\Cp{X_i}{G}$ by Lemma \ref{9.3}. It now
follows that $\Cp{X}{G}$ has precisely $n^m$ connected components.
Since each of these components is a translate of $c(\Cp{X}{G})$,
this finishes the proof of the ``only if'' part.

To prove the ``if'' part, assume that $c(\Cp{X}{G})$ has index $n^m$
in $\Cp{X}{G}$. If $X$ has finitely many connected components, then
the (finite) number of connected components of $X$ must be equal to
$m$. (This follows from the ``only if'' part of our proof.) So it
remains only to show that $X$ has finitely many connected
components. Assume the contrary. Then one can find pairwise disjoint
non-empty clopen sets $U_1,U_2,\ldots,U_{m+1}$ such that
$X=\bigcup_{i=1}^{m+1}U_i$. Observe that for $i=1,\ldots,m+1$, the
space $\Cp{U_i}{G}$ contains at least $n$-many pairwise disjoint
non-empty clopen subsets (namely, the sets $\Cp{U_i}{C_j}$ for
$j=1,\ldots, n$). Since $\Cp{X}{G}\cong\prod_{i=1}^{m+1}\Cp{U_i}{G}$
by Lemma \ref{9.3}, we conclude that $\Cp{X}{G}$ contains at least
$n^{m+1}$-many pairwise disjoint non-empty clopen subsets. Since
$n^{m+1}>n^m$, we obtain a contradiction.
\end{proof}

\begin{corollary}
\label{cor:souvislost2} Let $G$ be a topological group such that
$c_0(G)=c(G)$ is a proper subgroup of $G$ of finite index. Then
$G$-equivalence preserves the finite number of connected components.
\end{corollary}

The following example shows that there are groups $G$ which satisfy
neither the conditions of Proposition \ref{souvislost}, nor those of
Proposition \ref{souvislost2}, but for which $G$-equivalence
nevertheless preserves connectedness.

\begin{example}
\label{ex:ex} \rm{ {\em The group
$G=\mathbb{T}\times\prod_{i\in\mathbb{N}}\mathbb{Z}(i)$
has infinitely many elements of each order $n\ge 2$, and
$c_0(G)=c(G)~\mathbb{T}\times\{0\}$ has infinite index in $G$, yet
$G$-equivalence preserves the finite number of connected
components\/}. Indeed, since $\mathbb{T}$ is pathwise connected and
$\prod_{i\in\mathbb{N}}\mathbb{Z}(i)$ is hereditarily
disconnected, it follows from Theorem \ref{thm:4.4}(i)
 that $G$-equivalence implies $\mathbb{T}$-equivalence, and
$\mathbb{T}$-equivalence preserves the finite number of connected components by
Corollary \ref{cor:souvislost}. }
\end{example}

The last example in this section shows that the class of all groups
$G$ for which $G$-equivalence preserves
the finite number of connected components
 is not closed
under taking finite products,
even though this class is closed under taking finite powers by
Corollary
 \ref{connectedness:and:finite:powers}.
\begin{example}
\rm{ {\em If $G=\mathbb{T}\times\mathbb{Z}(2)^\omega$ and
$H=\mathbb{T}^\omega\times\mathbb{Z}(2)$, then both $G$-equivalence
and $H$-equivalence preserve the finite number of connected
components, but $(G\times H)$-equivalence does not\/}.
It follows from Theorem \ref{thm:4.4}(i) that $G$-equivalence
implies $\mathbb{T}$-equivalence and $H$-equivalence implies
$\mathbb{Z}(2)$-equivalence. Since both $\mathbb{T}$-equivalence
and $\mathbb{Z}(2)$-equivalence preserve the finite number of
connected components by Corollary \ref{cor:souvislost}, we conclude
that so do $G$-equivalence and $H$-equivalence.
On the other hand, since $G\times
H\cong(\mathbb{T}\times\mathbb{Z}(2))^\omega$, it follows from
Corollary \ref{notcon} that $(G\times H)$-equivalence does not
preserve the finite number of connected components.}
\end{example}

\section{$\mathscr{G}$-equivalence and its relation with
$G$-equivalence}
\label{scr:G:eq}

Definition \ref{reflexe},
Proposition \ref{existence:of:reflections}
and Proposition \ref{double:reflection} are well-known in category theory.
We include the proofs of these propositions only for the reader's convenience.

\begin{definition}\label{reflexe}
{\rm
\begin{itemize}
\item[(i)]
For a class $\mathscr{G}$ of topological groups
we denote by $\overline{\mathscr{G}}$ the smallest (with respect to
inclusion) class of topological
groups containing $\mathscr{G}$ which is closed under taking arbitrary
products and subgroups.

\item[(ii)]
Given a topological group $H$ and a class  $\mathscr{G}$ of
topological groups, we will say that $r_{H,\mathscr{G}}(H)\in\overline{\mathscr{G}}$
is a {\em reflection of $H$ in $\mathscr{G}$} provided that
there exists a continuous homomorphism
$r_{H,\mathscr{G}}:H\rightarrow
r_{H,\mathscr{G}}(H)$ (called the {\em reflection homomorphism})
satisfying the following condition:
For every $G\in\overline{\mathscr{G}}$ and each continuous homomorphism
$h: H\to G$ one can find a continuous homomorphism $g:r_{H,\mathscr{G}}(H)\rightarrow G$
such that $h=g\circ r_{H,\mathscr{G}}$.
\end{itemize}
}
\end{definition}

\begin{proposition}
\label{existence:of:reflections}
For every topological group $H$ and each class
$\mathscr{G}$ of topological groups the reflection $r_{H,\mathscr{G}}(H)$
of $H$ in $\mathscr{G}$ exists and is unique up to a
topological isomorphism.
\end{proposition}
\begin{proof}
There exists an indexed set
$\{(G_s,h_s): s\in S\}$ such that:
\begin{itemize}
\item[(a)] for each $s\in S$, $G_s\in \overline{\mathscr{G}}$ and
$h_s: H\to G_s$ is a continuous surjective homomorphism;
\item[(b)] if $G\in\overline{\mathscr{G}}$ and $h: H\to G$ is a continuous
homomorphism, then there exist $t\in S$, a subgroup $G_t'$ of $G$ and a topological isomorphism $i_t: G_t\to G_t'$ such that
$h=i_t\circ h_t$.
\end{itemize}
The diagonal product $r_{H,\mathscr{G}}=\triangle_{s\in S}
h_s:H\rightarrow\prod_{s\in S}G_s$ of the family $\{h_s:s\in S\}$ is
a continuous group homomorphism. Clearly, $r_{H,\mathscr{G}}(H)\in
\overline{\mathscr{G}}$. Let $G\in\overline{\mathscr{G}}$, and
let $h:
H\to G$ be a continuous homomorphism. Let $t\in S$ and $i_t$ be as
in the conclusion of item (b), and let $\pi_t: \prod_{s\in S}G_s\to
G_t$ be the projection on $t$'s coordinate. Then $g=i_t\circ
\pi_t\restriction_{r_{H,\mathscr{G}}(H)}:r_{H,\mathscr{G}}(H)\to
G_t'$ is a continuous group homomorphism such that $g\circ
r_{H,\mathscr{G}}=i_t\circ
\pi_t\restriction_{r_{H,\mathscr{G}}(H)}\circ\ r_{H,\mathscr{G}}
=i_t\circ h_t=h$. This proves the existence of
$r_{H,\mathscr{G}}(H)$.

To show its uniqueness, assume that $r_0: H\to r_0(H)\in\overline{\mathscr{G}}$
and
$r_1: H\to r_1(H)\in\overline{\mathscr{G}}$
are continuous homomorphisms such that, for every $i\in\{0,1\}$,
each $G\in\overline{\mathscr{G}}$ and every continuous homomorphism
$h: H\to G$ one can find a continuous homomorphism $g_{i,h}:r_i(H)\rightarrow G$
such that $h=g_{i,h}\circ r_i$.
In particular, $r_{1-i}=g_{i, r_{1-i}}\circ r_i$ for $i\in \{0,1\}$.
Fix $i\in \{0,1\}$. Then
$r_i=g_{1-i, r_{i}}\circ r_{1-i}=g_{1-i, r_{i}}\circ g_{i, r_{1-i}}\circ r_i$,
which yields that $g_{1-i, r_{i}}\circ g_{i, r_{1-i}}$ is the identity map on
$r_i(H)$.
Therefore, $g_{1, r_{0}}: r_1(H)\to r_0(H)$ is the inverse map of the map
 $g_{0, r_{1}}: r_0(H)\to r_1(H)$. Hence, $r_0(H)\cong r_1(H)$.
\end{proof}
\begin{proposition}
\label{double:reflection}
Suppose that $H$ is a topological group and $\mathscr{G}$, $\mathscr{G}'$ are
classes of topological groups. If $\mathscr{G}'\subseteq \mathscr{G}$, then
$r_{H,\mathscr{G}'}(H)\cong r_{\mathscr{G}', r_{H,\mathscr{G}}(H)}(r_{H,\mathscr{G}}(H))$.
\end{proposition}
\begin{proof}
Let $K=r_{H,\mathscr{G}}(H)$, $K'=r_{K,\mathscr{G}'}(K)$ and $H'=r_{H,\mathscr{G}'}(H)$.
We need to prove that
$H'\cong K'$.

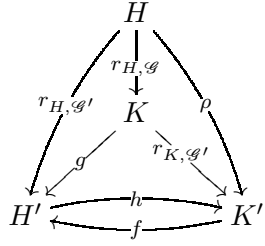
\begin{center} \medskip\hspace{1em}
\xymatrix{
&H\ar[d]|>>>>{r_{H,\mathscr{G}}}\ar@/_/[ldd]|{r_{H,\mathscr{G}'}}\ar@/^/[rdd]|{\rho}&\\
&K\ar[ld]|{g}\ar[rd]| <<<<< {r_{K,\mathscr{G}'}} 
&\\
H'\ar@/^/[rr]|h & &K'\ar@/^/[ll]|f}
\end{center}

Since
$H'\in\overline{\mathscr{G}'}\subseteq \overline{\mathscr{G}}$
and $K=r_{H,\mathscr{G}}(H)$, there exists a continuous homomorphism
$g: K\to H'$ such that $r_{H,\mathscr{G}'}=g\circ r_{H,\mathscr{G}}$.
Since $H'\in\overline{\mathscr{G}'}$ and
$K'=r_{K,\mathscr{G}'}(K)$, there exists a continuous homomorphism
$f:K'\to H'$ such that
$g=f\circ r_{K,\mathscr{G}'}$.
Define $\rho=r_{K,\mathscr{G}'}\circ r_{H,\mathscr{G}}: H\to K'$.
Then
$r_{H,\mathscr{G}'}=g\circ r_{H,\mathscr{G}}=f\circ r_{K,\mathscr{G}'}\circ r_{H,\mathscr{G}}=
f\circ \rho$.
Since $K'\in\overline{\mathscr{G}'}$ and $H'=r_{H,\mathscr{G}'}(H)$,
there exists a continuous homomorphism
$h: H'\to K'$ such that $\rho=h\circ r_{H,\mathscr{G}'}=h\circ f\circ \rho$.
Thus, $h\circ f$ is the identity map on $\rho(H)=K'$.
Finally,
$r_{H,\mathscr{G}'}=f\circ \rho=f\circ h\circ r_{H,\mathscr{G}'}$, which yields
that $f\circ h$ is the identity map on $r_{H,\mathscr{G}'}=H'$.
Therefore,
$h=f^{-1}$, and so $H'\cong K'$.
\end{proof}

\begin{definition}
\label{definition:of:free:object}
{\rm
For a space $X$ and a class $\mathscr{G}$ of topological groups
we define $F_\mathscr{G}(X)=r_{F(X),\mathscr{G}}(F(X))$
and call
$F_\mathscr{G}(X)$
the {\em free object over $X$ in
$\mathscr{G}$}.
}
\end{definition}

Since the free topological group $F(X)$ of $X$ is unique up to a
topological isomorphism, from Definition
\ref{definition:of:free:object} and Proposition
\ref{existence:of:reflections} we obtain the following

\begin{proposition}\label{uniqueA(X)}
For every space $X$ and each class $\mathscr{G}$ of topological
groups, the free object $F_\mathscr{G}(X)$ over $X$ in
$\mathscr{G}$ exists and is unique up
to a topological isomorphism.
\end{proposition}

One also has the ``usual properties'' of the free object:
\begin{proposition}
\label{freeobject}
For a space $X$ and a class $\mathscr{G}$ of topological groups,
there exists a continuous mapping
$\varphi_{X,\mathscr{G}}:X\rightarrow F_\mathscr{G}(X)$ satisfying the
following conditions:
\begin{enumerate}
\item[(i)]$\varphi_{X,\mathscr{G}}(X)$ algebraically generates
$F_\mathscr{G}(X)$;
\item[(ii)] for every $G\in\overline{\mathscr{G}}$ and each continuous
map $f:X\rightarrow G$ there exists a (unique) continuous
homomorphism $g:F_\mathscr{G}(X)\rightarrow G$ such that
$f=g\circ\varphi_{X,\mathscr{G}}$.
\end{enumerate}
\end{proposition}
\begin{proof}
Define
$\varphi_{X,\mathscr{G}}=r_{F(X),\mathscr{G}}\restriction_{X}$.
Since $X$ algebraically generates $F(X)$ and $r_{F(X),\mathscr{G}}$
is a homomorphism, we get (i). To prove (ii), assume that
$G\in\overline{\mathscr{G}}$ and $f:X\rightarrow G$ is a continuous
map. Let $h: F(X)\to G$ be the unique continuous homomorphism
extending $f$. Since $G\in\overline{\mathscr{G}}$ and
$F_\mathscr{G}(X)=r_{F(X),\mathscr{G}}(F(X))$, there exists a
continuous homomorphism $g: F_\mathscr{G}(X)\to G$ such that
$h=g\circ r_{F(X),\mathscr{G}}$. Now $f=h\restriction_X=g\circ
r_{F(X),\mathscr{G}}\restriction_X=g\circ \varphi_{X,\mathscr{G}}$.
The uniqueness of $g$ follows from (i).
\end{proof}

When
$\mathscr{G}$ forms a (wide) variety of topological groups,
the free object over $X$
in $\mathscr{G}$
 was
defined and investigated by S. Morris in \cite{Morris}. Comfort and
van Mill have generalized this concept in \cite{Comfort}.
In fact, when a space $X$ admits a homeomorphic embedding into a Cartesian product
of a family of members of $\mathscr{G}$, $F_\mathscr{G}(X)$ coincides with the free
topological group in the class $\mathscr{G}$ defined in \cite{Comfort}.

From Definition \ref{definition:of:free:object} and Proposition
\ref{double:reflection} one immediately gets the following
\begin{proposition}\label{VAR}
Let $X$ be a topological space and $\mathscr{G}$, $\mathscr{G}'$
two classes of topological groups.
If $\mathscr{G}'\subseteq\mathscr{G}$, then $F_{\mathscr{G}'}(X)\cong
r_{F_\mathscr{G}(X),\mathscr{G}'}(F_\mathscr{G}(X))$.
\end{proposition}

\begin{definition}
\label{definition:scr:G:equivalence}
{\rm
Let $\mathscr{G}$ be a class of topological groups.
We say that spaces $X$ and $Y$ are {\em $\mathscr{G}$-equivalent}
(and we write $X\eq{\mathscr{G}}Y$)
provided that $F_\mathscr{G}(X)\cong F_\mathscr{G}(Y)$.
}
\end{definition}

When $\mathscr{G}$ is the class of all topological groups, one
obviously has $F_\mathscr{G}(X)\cong F(X)$, and so
$\mathscr{G}$-equivalence in this case coincides with the classical
$M$-equivalence of Markov. Similarly, when $\mathscr{A}$ is the
class of all Abelian topological groups, then $F_\mathscr{A}(X)$
coincides with the free Abelian group in the sense of Markov, and so
$\mathscr{A}$-equivalence in this case coincides with the classical
$A$-equivalence of Markov.
Since $F_{\overline{\mathscr{G}}}(X)\cong F_{\mathscr{G}}(X)$,
$\overline{\mathscr{G}}$-equivalence is the same as $\mathscr{G}$-equivalence.

\begin{theorem}
If $\mathscr{G}'\subseteq\mathscr{G}$ are two classes of
topological groups, then $\mathscr{G}$-equivalence implies
$\mathscr{G}'$-equivalence.
\end{theorem}
\begin{proof}
Immediately follows from Definition
\ref{definition:scr:G:equivalence} and Proposition \ref{VAR}.
\end{proof}

Given topological groups $G$ and $H$, we denote by $\Chom_p(G,H)$
the subspace of $\Cp{G}{H}$ consisting of homomorphisms from $G$ to
$H$. If $H$ is Abelian, then $\Chom_p(G,H)$ is a topological group.

\begin{theorem}\label{X versus AX}
Assume that $\mathscr{G}$ is a class of
topological groups and $G\in\mathscr{G}$ is Abelian. Then
$\Cp{X}{G}\cong \Chom_p(F_\mathscr{G}(X),G)$ for every space $X$.
\end{theorem}
\begin{proof}
Let $\varphi_{X,\mathscr{G}}$ be as in Proposition \ref{freeobject}. For
every $f\in\Cp{X}{G}$ put $\phi(f)=\hat{f}$, where
$\hat{f}:F_\mathscr{G}(X)\to G$ is the
unique continuous homomorphism such that
$f=\hat{f}\circ\varphi_{X,\mathscr{G}}$.
Since $G$ is Abelian,
$\phi: \Cp{X}{G}\to \Chom_p(F_\mathscr{G}(X),G)$ is an
isomorphism.
We need
to show that $\phi$ is a homeomorphism as well.

If $x\in X$ and $V$ is an open
neighborhood of the identity in $G$, then
\begin{align*}
\phi(W(x,V))=&
\{\hat{f}:f\in\Cp{X}{G}, f(x)\in V\}
=
\{\hat{f}:f\in\Cp{X}{G},
\hat{f}(\varphi_{X,\mathscr{G}}(x))\in V\}\\
=&
\{\pi\in \Chom_p(F_\mathscr{G}(X),G):\pi(\varphi_{X,\mathscr{G}}(x))\in V\}
\end{align*}
is an open set in $\Chom_p(F_\mathscr{G}(X),G)$. Since the family
$\{W(x,V): x\in X, V$
is an open neighborhood of $e$ in $G\}$ forms a
subbase of open neighborhoods of the identity map of $\Cp{X}{G}$,
this proves that $\phi$ is an open map.

For $a\in
F_\mathscr{G}(X)$ and an open neighborhood $O$ of the identity in
$G$, define
$O_a=\{\pi\in \Chom_p(F_\mathscr{G}(X),G): \pi(a)\in O\}$.
Let us show that the set $\phi^{-1}(O_a)$ is open in $\Cp{X}{G}$.
Since
$\varphi_{X,\mathscr{G}}(X)$ generates $F_\mathscr{G}(X)$, there exist
$n\in\mathbb{N}\setminus\{0\}$, $x_1,\ldots,x_n\in X$ and $z_1,\ldots,z_n\in\mathbb{Z}$ such that
$a=\prod_{i=1}^n\varphi_{X,\mathscr{G}}(x_i)^{z_i}$.
Choose $f_0\in \phi^{-1}(O_a)$ arbitrarily.
Then $\phi(f_0)\in O_a$, and therefore
$$
\phi(f_0)(a)=\hat{f_0}(a)=\hat{f_0}\left(\prod_{i=1}^n\varphi_{X,\mathscr{G}}(x_i)^{z_i}\right)
=
\prod_{i=1}^n\left(\hat{f_0}(\varphi_{X,\mathscr{G}}(x_i))\right)^{z_i}
=
\prod_{i=1}^n f_0(x_i)^{z_i}
\in O.
$$
Using the continuity of group operations in $G$, for each $i=1,\dots,n$ we can choose an open
neighborhood $U_i$ of $f_0(x_i)$ in $G$ such
that $\prod_{i=1}^n U_i^{z_i}\subseteq O$.
Now
$W=\bigcap_{i=1}^n W(x_i,U_i)$
is an open subset of
$\Cp{X}{G}$ satisfying $f_0\in W\subseteq \phi^{-1}(O_a)$. This
proves that $\phi^{-1}(O_a)$ is open in $\Cp{X}{G}$. Since the
family $\{O_a:a\in F_\mathscr{G}(X), O$ is an open neighborhood of
$e$ in $G\}$ is a subbase of open neighborhoods of the identity in
$\Chom_p(F_\mathscr{G}(X),G)$, we conclude that $\phi$ is
continuous.
\end{proof}
Let us note that the condition ``$G$ is Abelian'' in the above theorem
cannot be omitted. Indeed, if $G$ is not Abelian,
$\Chom_p(F_\mathscr{G}(X),G)$ need not be a group,
 while
$\Cp{X}{G}$ is {\em always\/} a group.

\begin{corollary}\label{A_G implikuje G}
If $\mathscr{G}$ is a class of topological groups and
$G\in\mathscr{G}$ is Abelian, then $\mathscr{G}$-equivalence implies
$G$-equivalence.
\end{corollary}
\begin{proof}
If $X$ and $Y$ are $\mathscr{G}$-equivalent, then
$F_{\mathscr{G}}(X)\cong F_{\mathscr{G}}(Y)$, and so
$$
\Cp{X}{G}\cong \Chom_p(F_\mathscr{G}(X),G)\cong \Chom_p(F_\mathscr{G}(Y),G)\cong \Cp{Y}{G}
$$
by Theorem \ref{X versus AX}. Thus $X$ and $Y$ are $G$-equivalent.
\end{proof}

This corollary allows us to distinguish between $\{G\}$-equivalence
and $G$-equivalence for many Abelian topological groups $G$, as the
following example demonstrates.
\begin{example}\label{G:weaker:then:scrG}
Let $H$ be an Abelian topological group such that $H$-equivalence
preserves pseudocompactness within the class of $H$-regular spaces.
Let $G=H^\omega$. Then $\{G\}$-equivalence implies $G$-equivalence,
while $G$-equivalence does not imply $\{G\}$-equivalence. \rm{
Indeed, the first statement follows directly from Corollary \ref{A_G
implikuje G}. To prove the second statement, note that
$\overline{\{G\}}=\overline{\{H\}}$, so $\{G\}$-equivalence
coincides with $\{H\}$-equivalence. Hence $\{G\}$-equivalence
implies $H$-equivalence by Corollary \ref{A_G implikuje G}. Since
$H$-equivalence preserves pseudocompactness within the class of
$H$-regular spaces, so does $\{G\}$-equivalence. On the other hand,
it follows from Corollary \ref{Gomegadoesntpreserve} that
$G$-equivalence does not preserve pseudocompactness within the class
of $H$-regular spaces. Thus,
$G$-equivalence does not imply $\{G\}$-equivalence. }
\end{example}

\begin{lemma}\label{X closed in F(X)}
Let $\mathscr{G}$ be a class of topological groups and $X$ a
topological space such that $X$ is $G^\star$-regular for some
$G\in\mathscr{G}$.  Then $\varphi_{X,\mathscr{G}}:X\to
F_\mathscr{G}(X)$ is a homeomorphic embedding such that
$\varphi_{X,\mathscr{G}}(X)$  is closed in $F_\mathscr{G}(X)$.
\end{lemma}
\begin{proof}
Since $X$ is $G^\star$-regular
for some $G\in\mathscr{G}$,
 $\varphi_{X,\mathscr{G}}$  is a homeomorphic embedding. We
identify $X$ with $\varphi_{X,\mathscr{G}}(X)$. Fix $a\in
F_\mathscr{G}(X)\setminus X$. We have to prove that $a$ is not in
the closure of $X$. Since $X$ algebraically generates
$F_\mathscr{G}(X)$, there exist $n\in\mathbb{N}$, $x_1,\ldots,x_n\in
X$ and $z_1,\ldots,z_n\in\mathbb{Z}$ such that
$a=\prod_{i=1}^nx_i^{z_i}$. Define $s=\sum_{i=1}^nz_i$ if
$\{x_1,\ldots,x_n\}\not=\emptyset$ and let $s=0$ otherwise. We
consider two cases.

{\it Case 1\/}. {\sl There exist $G\in\mathscr{G}$ and $g\in G$ such
that $g^s\neq g$\/}. In this case take a constant map $f\in
C_p(X,G)$ defined by $f(x)=g$ for all $x\in X$. This map extends to
a continuous homomorphism $\hat{f}:F_\mathscr{G}(X)\to G$. Clearly,
$\hat{f}(a)=g^s\neq g$ and
$\hat{f}(X)=f(X)\subseteq\{g\}$,
so $V=\hat{f}^{-1}(G\setminus\{g\})$ is an open neighborhood of $a$
disjoint from $X$.

{\it Case 2\/}. {\sl $g^s=g$ for each $G\in\mathscr{G}$ and every $g\in G$\/}.
Fix pairwise disjoint open subsets $U$, $U_1,\dots,U_n$ of $F_\mathscr{G}(X)$
such that $a\in U$,
and $x_i\in U_i$ for $i=1,\ldots,n$.
Take $G\in\mathscr{G}$ such that $X$ is $G^\star$-regular. Then
there exists $g\in G\setminus\{e\}$ and $f_i\in C_p(X,G)$  such that
$f_i(x_i)=g$ and
$f_i(X\setminus U_i)\subseteq\{e\}$
 for every
$i=1,\ldots,n$. The function $f=\prod_{i=1}^nf_i$ extends to a
continuous homomorphism $\hat{f}:F_\mathscr{G}(X)\to G$.
Since
$\hat{f}(a)=g^s=g\not=e$,
$V=U\cap\hat{f}^{-1}(G\setminus\{e\})$
is an open neighborhood of $a$ in $F_\mathscr{G}(X)$.
If $x\in X\setminus \bigcup_{i=1}^nU_i$, then $\hat{f}(x)=e$,
so $x\not\in V$.
If $x\in\bigcup_{i=1}^nU_i$, then $x\not\in U$, so again $x\not\in V$.
This proves that $V\cap X=\emptyset$.

In both cases we have found an open neighborhood $V$ of $a$ which does
not intersect $X$. Since $a$ was chosen arbitrarily, we conclude
that $X$ is closed in $F_\mathscr{G}(X)$.
\end{proof}
\begin{proposition}\label{sigmacomp}
Let $X$ be a space and $\mathscr{G}$ a class of topological
groups such that $X$ is $G^\star$-regular for some
$G\in\mathscr{G}$. Then $X$ is $\sigma$-compact (Lindel\"of $\Sigma$-space) if and only if
$F_\mathscr{G}(X)$ is $\sigma$-compact (Lindel\"of $\Sigma$-space, respectively).
\end{proposition}
\begin{proof}
We start with the ``if'' part. By Lemma \ref{X closed in F(X)},
$X$ is a closed subspace of a $\sigma$-compact space (Lindel\"of
$\Sigma$-space) $F_\mathscr{G}(X)$. Hence $X$ is $\sigma$-compact
(Lindel\"of $\Sigma$-space, respectively).  Let us prove the ``only
if'' part. Since $X$ algebraically generates $F_\mathscr{G}(X)$,
there exists a representation $F_\mathscr{G}(X)=\bigcup_{i=0}^\infty
F_i$ where each $F_i$ is a continuous image of some finite power
$X^{n_i}$ of $X$. Now it remains only to note that the class of
$\sigma$-compact spaces (Lindel\"of $\Sigma$-spaces) is closed under
taking countable unions, finite powers and continuous images.
\end{proof}

\section{$\mathbb{T}$-equivalence}
\label{T:section}

In the theory of topological groups the group $\mathbb{T}$
can be viewed as
a
counterpart to $\mathbb{R}$ in the theory of topological vector
spaces. In this section we derive some corollaries about
$\mathbb{T}$-equivalence from theorems that we have established so far,
and we compare them with the similar statements about
$l$-equivalence ($\mathbb{R}$-equivalence), in order to emphasize
that the properties of $\mathbb{T}$-equivalence are at least as
good as (and often even better than) those of $l$-equivalence.

\begin{definition}
\label{dagger:definition}
{\rm
\begin{itemize}
\item[(i)]
We denote by $\mathscr{P}$ the class of all precompact Abelian
groups.
\item[(ii)]
For every topological group $G$, let
$G^\dagger=\Chom_p(G,\mathbb{T})$.
\end{itemize}
}
\end{definition}

The following ``precompact duality theorem'' is an immediate corollary
of the well-known results of Comfort and Ross
\cite[Theorems 1.2(c) and 1.3]{CR}, as well as
the particular case of a much more general result of Menini and Orsatti
(combine Propositions 2.8 and 3.9 with Theorem 3.11 in \cite{Menini}).
Additional information related to this theorem
can also be found in \cite{Racek} and
\cite{Hernandez}.

\begin{theorem}\label{duality}
$G\cong (G^\dagger)^\dagger$ for each $G\in\mathscr{P}$.
\end{theorem}

\begin{theorem}\label{dualT}
Let $X$ be a space. Then:
\begin{itemize}
\item[(i)]
 $\Cp{X}{\mathbb{T}}\cong
\Chom_p(F_\mathscr{P}(X),\mathbb{T})= F_{\mathscr{P}}(X)^\dagger$,
and
\item[(ii)]
$F_\mathscr{P}(X)\cong \Chom_p(\Cp{X}{\mathbb{T}},\mathbb{T})$.
\end{itemize}
\end{theorem}
\begin{proof}
(i)
Since
$\overline{\mathscr{P}}=\mathscr{P}$, from Definitions
\ref{reflexe} and \ref{definition:of:free:object}
it follows that
$F_\mathscr{P}(X)\in\mathscr{P}$. Since $\mathbb{T}\in\mathscr{P}$,
$\Cp{X}{\mathbb{T}}\cong
\Chom_p(F_\mathscr{P}(X),\mathbb{T})= F_{\mathscr{P}}(X)^\dagger$
by Theorem \ref{X versus AX}
and Definition \ref{dagger:definition}(ii).

(ii)
Applying
Theorem
\ref{duality}, item (i) and Definition \ref{dagger:definition}(ii),
we obtain
$F_\mathscr{P}(X)\cong
\left(F_\mathscr{P}(X)^\dagger\right)^\dagger
\cong
\Cp{X}{\mathbb{T}}^\dagger
=
\Chom_p(\Cp{X}{\mathbb{T}},\mathbb{T})$.
\end{proof}

It is well-known that $\Cp{X}{\mathbb{R}}\cong
\Chom_p(L_p(X),\mathbb{R})$ where $L_p(X)\cong
\Chom_p(\Cp{X}{\mathbb{R}},\mathbb{R})$; see \cite{A}. Consequently,
$X\eq{\mathbb{R}}Y$ if and only if $L_p(X)\cong L_p(Y)$.
Our next corollary establishes
a counterpart to this theorem for
$\mathbb{T}$.

\begin{corollary}\label{PjeT}
$\mathbb{T}$-equivalence coincides with $\mathscr{P}$-equivalence.
\end{corollary}
\begin{proof}
Let $X$ and $Y$ be spaces.

Assume that $X\eq{\mathbb{T}}Y$. Then
$\Cp{X}{\mathbb{T}}\cong\Cp{Y}{\mathbb{T}}$, and so
$F_{\mathscr{P}}(X)^\dagger\cong F_{\mathscr{P}}(Y)^\dagger$
by Theorem \ref{dualT}(i). Applying
Theorem \ref{duality}, we obtain
$F_{\mathscr{P}}(X)\cong (F_{\mathscr{P}}(X)^\dagger)^\dagger\cong (F_{\mathscr{P}}(Y)^\dagger)^\dagger\cong F_{\mathscr{P}}(Y)$.
Therefore, $X\eq{\mathscr{P}}Y$.

Assume now that $X\eq{\mathscr{P}}Y$. Then $F_{\mathscr{P}}(X)\cong F_{\mathscr{P}}(Y)$, and so
$F_{\mathscr{P}}(X)^\dagger\cong F_{\mathscr{P}}(Y)^\dagger$. Applying Theorem \ref{dualT}(i)
once again, we conclude that $\Cp{X}{\mathbb{T}}\cong\Cp{Y}{\mathbb{T}}$.
This proves $X\eq{\mathbb{T}}Y$.
\end{proof}

In general, $G$-equivalence is weaker then $\mathscr{G}$-equivalence
for an Abelian $G\in\mathscr{G}$; see Example
\ref{G:weaker:then:scrG}.
 However,  Corollary \ref{PjeT} shows this is not
the case for the class $\mathscr{P}$ of all precompact Abelian
groups and $\mathbb{T}\in \mathscr{P}$.

Recall that the compact group $\widehat{F_\mathscr{P}(X)}$ is called
the {\em free compact Abelian group\/} of a space $X$; see
\cite{Hofmann}. Our Corollary \ref{PjeT} should be compared with the
following result from \cite{Hofmann}: two spaces $X, Y$ generate the
same free compact Abelian group if and only if $C(X,\mathbb{T})$ and
$C(Y,\mathbb{T})$ are algebraically isomorphic.

\begin{corollary}
\label{T:implies:any:precompact}
$\mathbb{T}$-equivalence implies $G$-equivalence
for  every
precompact Abelian group $G$.
\end{corollary}

\begin{proof}
Combine Corollaries \ref{PjeT} and \ref{A_G implikuje G}.
\end{proof}

\begin{theorem}
\label{dual:properties:for:T}
Let $X$ be a space.
\begin{itemize}
\item[(i)]
$X$ is pseudocompact if and only if $\Cp{X}{\mathbb{T}}$ is
\IANS.
\item[(ii)]
$l^*(X)=t(\Cp{X}{\mathbb{T}})$.
\item[(iii)]
$X$ is compact if and only if $\Cp{X}{\mathbb{T}}$ is a
\IANS\ group of countable tightness.
\item[(iv)]
$X$ is compact metrizable if and only if $\Cp{X}{\mathbb{T}}$ is a
\IANS\ group with a countable network.
\item[(iv)]
$X$ is totally disconnected if and only if
$tor(\Cp{X}{\mathbb{T}})$ is dense in $\Cp{X}{\mathbb{T}}$.
\item[(vi)]
For a given integer $n\in\mathbb{N}\setminus\{0\}$, the space $X$ has
precisely $n$
connected components
if and only
if for every (equivalently, for some) prime number $p$ the group $C(X,\mathbb{T})$ has
exactly $p^n-1$ elements of order $p$.
\end{itemize}
\end{theorem}
\begin{proof}
Since $\mathbb{T}$ is pathwise connected, it follows from
Proposition \ref{Gregualrity} that every space $X$ is $\mathbb{T}^{\star\star}$-regular
(and thus, both $\mathbb{T}^{\star}$-regular and
$\mathbb{T}$-regular).
Since $\mathbb{T}$ is
a separable metric
 NSS
group,
 item (i) follows from Theorem \ref{main},
item (ii) follows from Corollary
 \ref{char:of:compactness} and item (iii) follows from
Proposition \ref{char:of:compact:metrizable}.
Item (iv) follows from Proposition \ref{tightness}, item (v) follows
from Theorem \ref{totdiscon}, and item (vi) follows from Proposition
\ref{souvislost}.
\end{proof}

\begin{theorem}
\label{properties:preserved:by:T}
$\mathbb{T}$-equivalence preserves the following properties:
\begin{itemize}
\item[(i)] pseudocompactness;
\item[(ii)] the cardinal invariant $l^*$ (defined in the beginning of Section \ref{seccomp});
\item[(iii)] property of being a Lindel\"of $\Sigma$-space;
\item[(iv)] $\sigma$-compactness;
\item[(v)] compactness;
\item[(vi)] the property of being compact metrizable;
\item[(vii)] the (finite) number of connected components;
\item[(viii)] connectedness;
\item[(ix)] total disconnectedness.
\end{itemize}
\end{theorem}
\begin{proof}
(i) and (ii) follow from items (i) and (iv) of Theorem \ref{dual:properties:for:T}, respectively.

(iii) and (iv) follow from Proposition \ref{sigmacomp} and Corollary
\ref{PjeT}.

(v) follows from items (i) and (iv).

(vi) follows from Theorem \ref{dual:properties:for:T}(iii).

(vii) follows from Theorem \ref{dual:properties:for:T}(vi).

(viii) follows from (vii).

(ix) follows from Theorem \ref{dual:properties:for:T}(v).
\end{proof}

From Corollary
\ref{A_G implikuje G} and Theorem \ref{properties:preserved:by:T} we immediately get:
\begin{corollary}
\label{scrG:equivalance:when:scrG:contains:T}
Let $\mathscr{G}$ be a
class
 of groups such that $\mathbb{T}\in\mathscr{G}$. Then
$\mathscr{G}$-equivalence preserves properties (i)--(ix) listed in
Theorem \ref{properties:preserved:by:T}.
\end{corollary}

For an infinite cardinal $\tau$, let $\mathscr{B}_\tau$ be the class
of topological groups $G$ such that for every open neighborhood $U$
of $e$ there exists a set $F\subseteq G$ with $G=F U$ and
$|F|<\tau$. Let $\mathscr{A}_\tau$ be the class consisting of
Abelian members of $\mathscr{B}_\tau$. Since $\mathbb{T}\in
\mathscr{A}_\tau\subseteq  \mathscr{B}_\tau$ for every infinite
cardinal $\tau$, from Corollary
\ref{scrG:equivalance:when:scrG:contains:T} we get
\begin{corollary}
\label{cor:11:9}
For every infinite cardinal $\tau$, both $\mathscr{A}_\tau$-equivalence
and $\mathscr{B}_\tau$-equivalence preserve each of the properties
(i)--(ix) listed in Theorem \ref{properties:preserved:by:T}.
\end{corollary}

To the best of our knowledge, the above corollary is
new
even for $\tau=\omega$. Similarly, the following particular case of
Corollary \ref{scrG:equivalance:when:scrG:contains:T} appears to be
new as well.
\begin{corollary}
\label{total:disconnectedness:M:A}
Total disconnectedness is preserved by $A$-equivalence and $M$-equivalence.
\end{corollary}

Note that $l$-equivalence preserves properties from items (i)--(vi)
of Theorem \ref{properties:preserved:by:T}. To the best of
the author's knowledge, it is not known whether $l$-equivalence
preserves total disconnectedness. On the other hand, it is known
that $l$-equivalence does not preserve properties from items (vii)
and (viii) of Theorem \ref{properties:preserved:by:T}. This allows us
to distinguish between $l$-equivalence and $\mathbb{T}$-equivalence.
\begin{proposition}
\label{R:does:not:imply:T}
$\mathbb{R}$-equivalence (that is, $l$-equivalence) does not imply $\mathbb{T}$-equivalence.
\end{proposition}
\begin{proof}
It is well-known that every connected metrizable space $X$ is
$l$-equivalent to the topological sum $X\oplus\{x\}$ of $X$ with a
singleton; see \cite{RP}. Hence $l$-equivalence does not preserve
connectedness.
Combining this with Theorem
\ref{properties:preserved:by:T}(viii),
we conclude that $l$-equivalence does not
imply $\mathbb{T}$-equivalence.
\end{proof}

\begin{remark}
{\rm
Our proof that $\mathbb{T}$-equivalence preserves
pseudocompactness necessarily differs from the existing proofs
that $\mathbb{R}$-equivalence preserves
pseudocompactness.
Indeed,
the original proof in \cite{Arkhangel} uses tools from
functional analysis (namely, barreled spaces) which are not available in the
case of $\mathbb{T}$. An alternative proof in \cite{Uspenskii}
is based on the following characterization:
a space $X$ is pseudocompact if and only if $\Cp{X}{\mathbb{R}}$ is
$\sigma$-precompact. This approach is not applicable in the
case of $\mathbb{T}$ since
$\Cp{X}{\mathbb{T}}$ is  precompact for {\em every\/} space $X$.
}
\end{remark}

\begin{remark}
{\rm We prove in our forthcoming paper \cite{SS} that the covering dimension
$\dim$ is preserved by $\mathbb{T}$-equivalence, so one can add this property
to the list of properties in Theorem \ref{properties:preserved:by:T}, as well
as in Corollaries \ref{scrG:equivalance:when:scrG:contains:T} and \ref{cor:11:9}.
}
\end{remark}

\section{Open questions}
\label{open:questions}

Besides general problems listed in Section \ref{introduction}, there
are many concrete questions that can be asked. In this section we
list only a small sample of those.

\begin{question}
\begin{itemize}
\item[(i)]
Can the assumption that $G$ is metric be omitted in Corollary
\ref{nademnou}?
\item[(ii)]
Can one replace ``$G^\star$-regularity'' by ``$G$-regularity'' in the assumption of  Corollary \ref{nademnou}?
\end{itemize}
\end{question}

\begin{question}
Is there an NSS group $G$ such that $G$-equivalence does
not preserve compactness within the class of $G$-regular  spaces?
\end{question}

\begin{question}
Is there a \IANS\ group $G$ such that $G$-equivalence does not preserve compactness
(or pseudocompactness) within the class of $G$-regular  spaces?
\end{question}

The following question arises in connection with Proposition
\ref{R:does:not:imply:T}:
\begin{question}
Does $\mathbb{T}$-equivalence imply $\mathbb{R}$-equivalence (that is, $l$-equivalence)?
\end{question}

\begin{question}
Let $\mathscr{R}=\{\mathbb{R}\}$ be the class consisting of a single
group $\mathbb{R}$ of reals. Do $\mathscr{R}$-equivalence and
$\mathbb{R}$-equivalence (that is, $l$-equivalence) coincide?
\end{question}

Note that a similar question for the torus $\mathbb{T}$ has a
positive answer. Indeed, $\overline{\{\mathbb{T}\}}=\mathscr{P}$,
and since $\mathbb{T}$-equivalence coincides with
$\mathscr{P}$-equivalence by Corollary \ref{PjeT}, it follows that
$\mathbb{T}$-equivalence coincides with $\mathscr{T}$-equivalence,
where $\mathscr{T}=\{\mathbb{T}\}$ is the class consisting of a
single group $\mathbb{T}$.
In connection with Corollary \ref{T:implies:any:precompact}, the
following question seems to be interesting:

\begin{question}
\begin{itemize}
\item[(i)] Does
$\mathbb{Q}^*$-equivalence
 imply $\mathbb{T}$-equivalence?
(Here $\mathbb{Q}^*$ denotes the Pontryagin dual of the discrete
group $\mathbb{Q}$ of rational numbers.)
\item[(ii)] For a prime number $p$, does $\mathbb{Z}_p$-equivalence imply
$\mathbb{T}$-equivalence within the class of
$\mathrm{ind}$-zero-dimensional Tychonoff spaces?
\end{itemize}
\end{question}

Let us recall another notion from the $C_p$-theory. Spaces $X$ and
$Y$ are called {\em $u$-equivalent\/} provided that the topological
groups $\Cp{X}{\mathbb{R}}$ and $\Cp{Y}{\mathbb{R}}$ considered with
their (two-sided) uniformities, are uniformly homeomorphic (that is,
there exists a bijection between $\Cp{X}{\mathbb{R}}$ and
$\Cp{Y}{\mathbb{R}}$ that preserves the uniform structures). This
motivates the following definition.

\begin{definition}
{\rm
For a given topological group $G$, let us say that spaces $X$ and $Y$ are {\em $G_u$-equivalent\/} ({\em $G_u'$-equivalent\/})
provided that topological groups   $\Cp{X}{G}$ and $\Cp{Y}{G}$ considered with their left uniformities
(two-sided uniformities, respectively), are uniformly homeomorphic.
}
\end{definition}

This definition leads to natural ``uniform analogues'' of our Problems \ref{problem:1}
 and \ref{problem:2}.
\begin{problem}
Given a topological group $G$, characterize topological properties of
$X$ in terms of the uniform properties of
$\Cp{X}{G}$.
\end{problem}
\begin{problem}
Given a topological group $G$, a class $\mathscr{C}$ of spaces and a
topological property $\mathscr{E}$, investigate when the property $\mathscr{E}$
is preserved
by $G_u$-equivalence ($G_u'$-equivalence) within the class $\mathscr{C}$.
\end{problem}

Remark \ref{compact-like:TAP:vs:NSS} motivates the following question:
\begin{question}
{\rm \cite{DSS}}
Is there a ZFC example of a countably compact Abelian \IANS\ group
that is not NSS?
\end{question}

\begin{question}
\label{question:G:regularity}
Is there a space $X$ and a topological group $G$ such that $X$ is
$G$-regular but not $G^\star$-regular?
\end{question}

\medskip
\noindent
{\bf Acknowledgement:\/} We would like to thank Dikran Dikranjan
for helpful discussions related to categorical notions in the area between \ref{reflexe} and
\ref{VAR},
and Oleg Okunev for supplying the reference \cite[Theorem 7]{R-quotient}.
We are indebted to G\'abor Luk\'acs for pointing out an incorrect bibliographic reference
related to Theorem \ref{duality} in the first version of this preprint.


\begin{thebibliography}{99}
\bibitem{A}
A.~V.~Arhangel'ski\u{\i}, {\em Topological function spaces\/},
Math.Appl., Vol.~78, Kluwer Academic, Dordrecht (1992).


\bibitem{Arkhangel}
A.~V.~Arhangel'ski\u{\i}, {\em On linear homeomorphisms
of function spaces\/}, Soviet. Math. Dokl. ~25 (1982), no. 3,
852-855. (Translated from Russian)

\bibitem{RP}
A.~V.~Arhangel'ski\u{\i},
{\em $C_p$-Theory\/}, in: Recent
Progress in General Topology, Elsevier Science Publishers B.V.
(1992).

\bibitem{R-quotient}
A.~V.~Arhangel'ski\u{\i},
{\em $R$-factor mappings of spaces with a countable base\/},
Dokl. Akad. Nauk SSSR ~287 (1986), no. 1, 14--17. (in Russian);
English
translation:
Soviet Math. Dokl.~33 (1986), no. 2, 302--305.

\bibitem{AT}
A.~V.~Arhangel'ski\u{\i}, M.~Tkachenko,
{\em Topological groups and related structures\/} (Atlantis Studies in Mathematics, 1), Atlantis Press, Paris; World Scientific Publishing Co. Pte. Ltd., Hackensack, NJ, 2008.

\bibitem{Comfort}W.~W.~Comfort, J.~van~Mill, {\em On the existence
of  free topological groups\/}, Topology Appl.~29 (1988),
245--269.

\bibitem{CR}W.~W.~Comfort, K.~A.~Ross, {\em Topologies induced by groups of characters\/}, Fund. Math.~55
(1964), 283--291.

\bibitem{Dikran} D.~N.~Dikranjan, I.~R.~Prodanov, L.~N.~Stoyanov,
{\em Topological groups: Characters, Dualities and Minimal
Group Topologies\/}, Pure and Applied Mathematics, Vol.~130,
Marcel Dekker Inc., New York-Basel, (1989).

\bibitem{DS}
D.~Dikranjan, D.~Shakhmatov,
{\em Selected topics from the structure theory of topological groups\/}, pp.
389-406 in: Open Problems in Topology II (E.~Pearl, Eds.) (North-Holland, Amsterdam, 2007).

\bibitem{DSS}
D.~Dikranjan, D.~Shakhmatov, J.~Sp\v{e}v\'{a}k,
{\em NSS and TAP properties in topological groups close to being compact\/},
preprint no.~arXiv:0909.2381 [math.GN].

\bibitem{DT}
X.~ Dom\'{\i}nguez, V.~Tarieladze,
{\em Metrizable TAP, HTAP and STAP groups\/},
preprint no.~arXiv:0909.1400 [math.GN].


\bibitem{Engelking}R.~Engelking, {\em General topology\/}, Heldermann,
(1989).

\bibitem{Graev}M.~I.~Graev, {\em Free topological groups\/}, Izv. Akad.
Nauk. SSSR Ser. Mat.~12~(3) (1948), 279-324 (in Russian); English
translation: Amer. Math. Soc. Transl.~(1) 8 (1962).

\bibitem{Guran}I.~I.~Guran, {\em Topological groups similar to Lindel\"{o}f groups\/}
(Russian), Dokl. Akad. Nauk SSSR~256 (1981), no. 6, 1305-1307.

\bibitem{Hernandez}S.~Hern\'{a}ndez, S.~Macario, {\em Dual properties in totally
bounded Abelian groups\/}, Arch. Math. (Basel)~80 (2003), no.
3, 271-283.

\bibitem{hern}S.~Hern\'{a}ndez, A.~M.~R\'{o}denas, {\em Automatic
continuity and representation of group homomorphisms defined between
groups of continuous functions\/}, Topology Appl.~154
(2007), 2089-2098.

\bibitem{Hofmann}K.~H.~Hofmann, S.~A.~Morris, {\em Free compact groups
I: Free compact Abelian groups\/}, Topology Appl.~23 (1986),
41-64.

\bibitem{Markov}A.~A.~Markov, {\em On free topological groups\/},
Izv. Akad. Nauk SSSR Ser. Mat.~9 (1) (1945), 3-64 (In Russian);
English translation: Amer. Math. Soc. Transl.~(1) 8 (1962).

\bibitem{Menini}C.~Menini, A.~Orsatti, {\em Dualities between
categories of topological modules\/}, Communications in Algebra~11~(1) (1983), 21-66.

\bibitem{Morris}S.~A.~Morris, {\em Varieties of topological groups: a survey\/},
 Colloq. Math.~46 (1982), no. 2, 147--165.

\bibitem{Pytkeev}E.~G.~Pytkeev, {\em On the tightness of spaces of
continuous functions\/}, Russian Math. Surveys~37: 1 (1982),
176-177. (Translated from Russian)

\bibitem{Racek}S.~U.~Raczkowski, F.~J.~Trigos-Arrieta, {\em Duality
of totally bounded Abelian groups\/}, Bol. Soc. Mat. Mexicana~(3)~ 7 (2001), no. 1, 1-12.

\bibitem{SS}
D.~Shakhmatov, J.~Sp\v{e}v\'{a}k,
{\em $\mathbb{T}$-equivalence preserves the covering dimension\/}, work in progress.

\bibitem{Sirota}
S.~M.~Sirota,
{\em A product of topological groups, and extremal disconnectedness\/} (in Russian),
Mat. Sb. (N.S.)~79~(121) (1969), 179--192.

\bibitem{Tkachuk}V.~V.~Tkachuk, {\em Duality with respect to the functor $C_p$ and cardinal invariants of the type of the Suslin
number\/} (Russian), Mat. Zametki~37 (1985), No.3, 441-451.

\bibitem{Uspenskii}V.~V.~Uspenski\u{\i}, {\em A characterization of
compactness in terms of the uniform structure in the space of
functions\/}, Uspekhi Matem. Nauk~37~(4) (1982), 183-184.
\end{thebibliography}
\end{document}